\documentclass[a4paper,10pt]{article}
\usepackage{amssymb}
\usepackage{latexsym}
\usepackage{amsmath}
\usepackage{amsthm}
\usepackage{amsfonts}
\usepackage{amssymb}
\usepackage{amsmath,amscd}
\usepackage{graphicx}

\newtheorem{Th}{Theorem}
\newtheorem{Prop}{Proposition}

\newtheorem{Lm}{Lemma}
\newtheorem{Lma}{Lemma}[section]

\newtheorem{Rm}{Remark}

\newcommand{\be}{\begin{equation}}
\newcommand{\ee}{\end{equation}}
\newcommand{\bes}{\begin{equation*}}
\newcommand{\ees}{\end{equation*}}

\newcommand{\R}{\mathbb{R}}
\newcommand{\N}{\mathbb{N}}
\newcommand{\C}{\mathbb{C}}

\newcommand\res{\mathop{\hbox{\vrule height 7pt width .5pt depth 0pt
\vrule height .5pt width 6pt depth 0pt}}\nolimits}

\def\theequation{\thesection.\arabic{equation}}
\def\theTh{\Roman{section}.\arabic{Th}}
\def\theProp{\Roman{section}.\arabic{Prop}}
\def\theCo{\Roman{section}.\arabic{Co}}

\def\theRm{\Roman{section}.\arabic{Rm}}

\newcommand{\reset}{\setcounter{equation}{0}\setcounter{Th}{0}\setcounter{Prop}{0}\setcounter{Co}{0}
\setcounter{Lm}{0}\setcounter{Rm}{0}}

\def\ti{\tilde}
\def\lf{\left}
\def\rg{\right}

\def\al{\alpha}
\def\la{\lambda}

\def\ep{\varepsilon}

\def\ds{\displaystyle}
\def\ov{\overline}
\def\Om{\Omega}
\def\om{\omega}
\def\p{\partial}

\def\res{\mathop{\hbox{\vrule height 7pt width .5pt 
depth 0pt\vrule height .5pt width 6pt depth 0pt}}\nolimits}

\begin{document}

\title{Lower Semi-Continuity of the Index in the Visosity Method for Minimal Surfaces.}

\author{ Tristan Rivi\`ere\footnote{Department of Mathematics, ETH Zentrum,
CH-8093 Z\"urich, Switzerland.}}
\maketitle

{\bf Abstract :}{\it The goal of the present work is twofold. First we prove the existence of an {\it Hilbert Manifold} structure on the space of  immersed oriented closed surfaces with three derivatives in $L^2$ in an arbitrary
sub-manifold $M^m$ of an euclidian space ${\R}^Q$. Second, using this Hilbert manifold structure, we prove a lower semi continuity property of the index for sequences of  conformal immersions, critical points to the viscous approximation of the area satisfying Struwe entropy estimate and bubble tree strongly converging in $W^{1,2}$ to a limiting minimal surface as the viscous parameter is going to zero.} 

\medskip

\noindent{\bf Math. Class. 49Q05, 53A10, 49Q10}

\section{Introduction}
Let $M^m$ be a smooth $m-$dimensional sub-manifold of an euclidian space ${\R}^Q$ and denote by $\pi_{M^m}$ the orthogonal projection onto $M^m$ defined in a neighborhood
of $M^m$.

 Let $\Sigma$ be an arbitrary closed oriented 2-dimensional manifold. We define the {\it Sobolev Space}
of maps between $\Sigma$ and $M^m$ with 3 derivatives in $L^2$ as follows
\[
W^{3,2}(\Sigma, M^m):=\lf\{   u\in W^{3,2}(\Sigma,{\R}^Q)\quad;\quad u(x)\in M^m\quad \forall\, x\in \Sigma  \rg\}
\]
where the Sobolev space of $W^{3,2}$ functions on $\Sigma$ is defined with respect to any arbitrary reference metric (they are all equivalent due to the compactness of $\Sigma$).

Since $ 3\times 2=6>2=\mbox{dim}(\Sigma)$ the space $W^{3,2}(\Sigma, M^m)$ inherits a natural {\it Hilbert Manifold} structure (see \cite{Riv-minimax}). Within this manifold
we are considering the open subset of {\it $W^{3,2}-$immersions}
\[
\mbox{Imm}^{3,2}(\Sigma, M^m):=\lf\{   \vec{\Phi}\in W^{3,2}(\Sigma, M^m)\quad;\quad d\vec{\Phi}\wedge d\vec{\Phi}\ne 0\quad \mbox{ in }\Sigma  \rg\}
\]
Since there is no ambiguity on the regularity we are choosing we shall simply omit the superscripts $3,2$ and denote $\mbox{Imm}(\Sigma, M^m)$ for $\mbox{Imm}^{3,2}(\Sigma, M^m)$. For a given oriented closed surface $\Sigma$ we denote by $b(\sigma)$ the sum of the first 3 first {\it Betti Numbers} of $\Sigma$ :
\[
b(\Sigma):=b_0(\Sigma)+b_1(\Sigma)+b_2(\Sigma)\quad.
\]
There are obviously finitely many classes modulo diffeomorphisms of surfaces $\Sigma$ such that $b(\Sigma)\le b$. We will work from now on with one fixed representative
in each of these classes. 

For any $b\in {\N}^\ast$ we denote by $\mbox{Imm}_b(M^m)$ the  {\it Hilbert Manifold} obtained by taking the disjoined union of the  {\it Hilbert Manifolds} of $W^{3,2}-$immersions of the finitely many surfaces such that $b(\Sigma)\le b$.

 Starting from $\mbox{Imm}_b(M^m)$ we are constructing an {\it Hilbert manifold} of  $W^{3,2}$ immersed surfaces in the following way. We are first marking each surface $\Sigma$ by fixing respectively 3, 1 or zero  distinct points on each component of $\Sigma$ of genus respectively 0, 1 and $>1$. We are then considering the quotients of  $\mbox{Imm}(\Sigma,M^m)$ by $\mbox{Diff}^{\,\ast}_+(\Sigma)$, the  positive $W^{3,2}-$diffeomorphisms of  $\Sigma$ preserving the points we have fixed and isotopic to the identity. Then we denote by 
 \[
{\mathfrak M}_b(M^m):=\bigsqcup_{ b(\Sigma)\le b}\mbox{Imm}(\Sigma,M^m)/\mbox{Diff}^{\,\ast}_+(\Sigma)\quad \mbox{ and }\quad{\mathfrak M}(M^m):=\bigsqcup_{ b(\Sigma)<+\infty}\mbox{Imm}(\Sigma,M^m)/\mbox{Diff}^{\,\ast}_+(\Sigma)
 \]
Our first main result is the following
\begin{Th}
\label{th-1}
For any $b\in {\N}\cup\{\infty\}$ there exists a {\it Hilbert Manifold} Structure on ${\mathfrak M}_b(M^m)$  such that the canonical projection
\[
\Pi \ : \ \mbox{Imm}_b(M^m)\ \longrightarrow\ {\mathfrak M}_b(M^m)
\]
is a smooth map. \hfill $\Box$
\end{Th}
Since $\mbox{Diff}^{\, \ast}_+(\Sigma^g)$ misses
to be a {\it Banach Lie Group} but is only a {\it topological Group}\footnote{On the space of $W^{3,2}-$diffeomorphisms the right multiplication is smooth but the left multiplication is not differentiable, the inverse mapping is not $C^1$, there is no canonical chart in the neighborhood of the identity, the exponential map is continuous but not $C^1$, it
is not locally surjective in a neighborhood of the identity, the Bracket operation in the Tangent space to the identity is not continuous....etc see a description of all these ``pathological behavior'' for instance in \cite{EbMa} or \cite{Mil}} with an {\it Hilbert Manifold} structure the existence of a differentiable Hilbert structure on the quotient $\mbox{Imm}_0(\Sigma^g,M^m)/\mbox{Diff}^{\, \ast}_+(\Sigma^g)$ is not the result of classical
consideration and deserves to be studied with care (Progresses in this direction are given in \cite{BiFi} for $W^{3,2}-$embeddings but we are not going to follow this approach and the one
we choose is more specific to surfaces but more precise too). 

We shall in fact make the previous theorem more precise and to that aim we introduce some notations.
Let $\Sigma$ be a closed oriented 2-dimensional manifold, and $\vec{\Phi}\in W_{imm}^{3,2}(\Sigma, M^m)$ and let
$g_{\vec{\Phi}}:=\vec{\Phi}^{\,\ast} g_{M^m}$. Denote by $\wedge^{1,0}\Sigma$ the canonical bundle of 1-0 forms over the Riemann surface issued from $(\Sigma,g_{\vec{\Phi}})$ and denote by
$P_{\vec{\Phi}}$ the $L^2$ projection orthogonal projection from $(\wedge^{(1,0)}\Sigma)^{\otimes^2}$ onto the space  of {\it holomorphic quadratic forms} $\mbox{Hol}_Q(\Sigma,g_{\vec{\Phi}})$ on $(\Sigma,g_{\vec{\Phi}})$  and by $P_{\vec{\Phi}}^\perp:=\mbox{Id}-P_{\vec{\Phi}}$. Define the linear map
\[
\begin{array}{rcl}
\ds {D}^\ast_{\vec{\Phi}}\ :\  T_{\vec{\Phi}}\mbox{Imm}(\Sigma,M^m) &\longrightarrow &\ds W^{2,2}\lf(   (\wedge^{(1,0)}\Sigma)^{\otimes^2} \rg)\\[5mm]
\ds \vec{w} &\longrightarrow & P_{\vec{\Phi}}^\perp\lf({\p}\vec{w}\ \dot{\otimes}\ {\p}\vec{\Phi}\rg)
\end{array}
\]
where, in local complex coordinates for $\vec{\Phi}$ we denote
\[
 {\p}\vec{w}\ \dot{\otimes}\ {\p}\vec{\Phi}:={\p}_{{z}}\vec{w}\, \cdot\, {\p}_{{z}}\vec{\Phi}\ \ d{z}\otimes d{z}
\]
We are now going to prove the following theorem 
\begin{Th}
\label{th-hilbert}
Let $\vec{\Phi}\in\mbox{Imm}(\Sigma,M^m)$, then there exists an open neighborhood ${\mathcal O}_{\vec{\Phi}}$ of $\vec{\Phi}$ in $\mbox{Imm}(\Sigma,M^m)$ invariant under the action of $\mbox{Diff}^{\, \ast}_+(\Sigma)$ and
two smooth maps on ${\mathcal O}_{\vec{\Phi}}$, equivariant under the action of $\mbox{Diff}^{\, \ast}_+(\Sigma)$,
\[
\lf\{
\begin{array}{rcl}
\ds{\vec{w}}_{\vec{\Phi}}\ :\ {\mathcal O}_{\vec{\Phi}} &\longrightarrow & \mbox{Ker}( {D}^\ast_{\vec{\Phi}})\subset T_{\vec{\Phi}}\mbox{Imm}(\Sigma,M^m) \\[5mm]
\ds{\Psi}_{\vec{\Phi}}\ :\ {\mathcal O}_{\vec{\Phi}} & \longrightarrow & \mbox{Diff}^{\, \ast}_+(\Sigma)
\end{array}
\right.
\]
satisfying 
\[
\forall\ \vec{\Xi}\in {\mathcal O}_{\vec{\Phi}}\quad\quad\vec{\Xi}\circ\Psi_{\vec{\Phi}}(\vec{\Xi})=\pi_{M^m}\lf(\vec{\Phi}+\vec{w}_{\vec{\Phi}}(\vec{\Xi})\rg)
\]
where $\pi_{M^m}$ is the orthogonal projection onto $M^m$ and ${\mathcal T}_{\vec{\Phi}}:=({\vec{w}}_{\vec{\Phi}},{\Psi}_{\vec{\Phi}})$ realizes a diffeomorphism from ${\mathcal O}_{\vec{\Phi}}$ onto  $U_{\vec{\Phi}}\times \mbox{Diff}^{\, \ast}_+(\Sigma)$ where $U_{\vec{\Phi}}$ is a neighborhood of  $0$ in $\mbox{Ker}( \ov{D}^\ast_{\vec{\Phi}})$. Moreover the map ${\mathcal T}_{\vec{\Phi}}$ satisfies the following equivariance property : $\forall \, \vec{\Xi}\in {\mathcal O}_{\vec{\Phi}}$ and forall $\Psi_0\in \mbox{Diff}^{\, \ast}_+(\Sigma)$ one has
\[
\Psi_{\vec{\Phi}}(\vec{\Xi}\circ{\Psi}_0)=\Psi_0^{-1}\circ\Psi_{\vec{\Phi}}(\vec{\Xi})\quad\mbox{ and }\quad\vec{w}_{\vec{\Phi}}(\vec{\Xi}\circ\Psi_0)=\vec{w}_{\vec{\Phi}}(\vec{\Xi})\quad.
\]
The space $\mbox{Imm}(\Sigma,M^m)/\mbox{Diff}^{\, \ast}_+(\Sigma)$ is Hausdorff and defines a Hilbert Manifold such that the projection map
\[
\mbox{Imm}(\Sigma,M^m)\longrightarrow{\mathfrak M}_{\Sigma}(M^m):=\mbox{Imm}(\Sigma,M^m)/\mbox{Diff}^{\, \ast}_+(\Sigma)
\]
defines a $\mbox{Diff}^{\, \ast}_+(\Sigma)-$bundle for which $({\mathcal T}_{\vec{\Phi}})_{\vec{\Phi}}$ represents a local trivialization. \hfill $\Box$\end{Th}

\begin{Rm}
\label{rm-coulomb}
The condition 
\[
D^\ast_{\vec{\Phi}}\vec{w}:=P_{\vec{\Phi}}^\perp\lf({\p}\vec{w}\ \dot{\otimes}\ {\p}\vec{\Phi}\rg)=0
\]
corresponds\footnote{Observe that ${\p}\vec{w}\ \dot{\otimes}\ {\p}\vec{\Phi}=0$ is the condition which, starting from a conformal immersion $\vec{\Phi}$, preserves the conformality of $\vec{\Phi}+t\vec{w}$
at the first order for the same riemman structure on $\Sigma$. Similarly, $d^\ast_Aa=0$ is the first order condition which, starting from a Coulomb gauge $d_A^\ast A=0$ preserve the Coulomb condition for $A+ta$ for the same covariant co-differentiation $d^\ast_A$. Moreover, it is well known that the conformality of an immersion $\vec{\Phi}$ can be interpreted as a Coulomb condition (see for instance \cite{Riv-park-city}) with respect to the action of the ``gauge group'' $\mbox{Diff}^{\, \ast}_+(\Sigma)$.}
 to the Coulomb slice condition
\[
d^\ast_Aa=0
\]
in Gauge Theory while studying the Hilbert Bundle structure of the quotient of $H^s$ connections by the Gauge group for $s>n/2$ and away from reducible connections (see \cite{FU}).
\hfill $\Box$
\end{Rm}

Once this {\it Hilbert Bundle} structure will be established we shall be considering the following application to the {\it viscosity method} for the area functional introduced by the author
in \cite{Riv-minmax}.

For any immersion $\vec{\Phi}\in \mbox{Imm}(\Sigma,M^m)$ we denote
\[
F(\vec{\Phi}):=\int_{\Sigma}\lf[1+|{\vec{\mathbb I}}_{\vec{\Phi}}|^2\rg]^2\ dvol_{g_{\vec{\Phi}}}
\]
where ${\vec{\mathbb I}}_{\vec{\Phi}}=\pi_{\vec{n}}(\nabla^Md\vec{\Phi})$ is the second fundamental form of the immersion $\vec{\Phi}$ in $M^m$. 

\medskip

Observe that 
\be
\label{I.001}
\forall\ b\in{\N} \quad\quad\exists\ C_{b}>0\quad\quad F(\vec{\Phi})<C_{b}\quad\Longrightarrow\quad\mbox{b}\,(\Sigma)\le b
\ee
This is a direct consequence of {\it Gauss Bonnet theorem} and {\it Cauchy Schwartz inequality}.

It is clear that $F(\vec{\Phi})$ only depends on the equivalence class $[\vec{\Phi}]$
of $\vec{\Phi}$ in ${\mathfrak M}_\Sigma(M^m)$. Since $F$ is a smooth functional on $\mbox{Imm}(\Sigma,M^m)$ (see \cite{Riv-minmax}) it descends to a smooth functional on 
${\mathfrak M}_\Sigma(M^m)$. We shall  prove the following theorem.
\begin{Prop}
\label{th-fredholm}
Let $[\vec{\Phi}]$ be a critical point of $F$ in ${\mathfrak M}(M^m)$. Then the second derivative of $F$ at $[\vec{\Phi}]$ defines a Fredholm and Elliptic operator .\hfill $\Box$
\end{Prop}
The {\it viscosity method} consists in studying the variations of the {\it area lagrangian}
\[
\mbox{Area}(\vec{\Phi})=\int_{\Sigma_{\vec{\Phi}}}\, dvol_{g_{\vec{\Phi}}}
\]
by considering relaxations of the form
\[
A^\sigma(\vec{\Phi}):=\mbox{Area}(\vec{\Phi})+\sigma^2\, F(\vec{\Phi})=\mbox{Area}(\vec{\Phi})+\sigma^2\, \int_{\Sigma_{\vec{\Phi}}}\lf[1+|{\vec{\mathbb I}}_{\vec{\Phi}}|^2\rg]^2\ dvol_{g_{\vec{\Phi}}}\quad.
\]
where $\sigma>0$. The work \cite{Riv-minmax}  has been devoted to the asymptotic analysis of sequences of critical points of $A^{\sigma_k}$, with uniformly bounded 
$A^{\sigma_k}$ energy and satisfying {\it Struwe's entropy condition}
\be
\label{entropy}
\sigma_k^2\, F(\vec{\Phi}_{k})= o\lf(\frac{1}{\log\sigma_k^{-1}}\rg)\quad\mbox{ as $\sigma_k$ goes to zero.}
\ee
It is proved in this two work that, modulo extraction of a subsequence, the immersions $\vec{\Phi}_{k}$ {\it varifold converges} towards a 2-dimensional {\it integer rectifiable stationary varifold} $\mathbf v_\infty$ of $M^m$ which is parametrized. In \cite{Riv-reg} and \cite{PiR} the regularity of the {\it parametrized integer rectifiable stationary varifold} is established. Hence we have the following theorem.
\begin{Th}
\label{th-visco}[\cite{Riv-minmax}, \cite{Riv-reg}, \cite{PiR}]
Let $\vec{\Phi}_k$ be a sequence of immersions of a closed surface $\Sigma$, critical points of $A^{\sigma_k}$ and such that
\[
\limsup_{k\rightarrow+\infty}A^{\sigma_k}(\vec{\Phi}_k)<+\infty\quad\mbox{and }\quad\sigma_k^2\int_{\Sigma^g}(1+|\vec{\mathbb I}_{\vec{\Phi}_k}|^2)^2\ dvol_{g_{\vec{\Phi}_k}}=o\lf( \frac{1}{\log\sigma_k^{-1}} \rg)
\]
Then there exists a subsequence $\vec{\Phi}_{k'}$ such that the corresponding associated varifold\footnote{The associated varifold $\mathbf v$ associated to an immersion $\vec{\Phi}$ of $\Sigma^g$ is given by
\[
\forall \ \phi\in C^0(G_2TM^m)\quad {\mathbf v}(\phi):=\int_{\Sigma^g} \phi(\vec{\Phi}_\ast T_x \Sigma^g)\ dvol_{\vec{\Phi}^{\,\ast} g_{M^m}}
\]
}of  ${\mathbf v}_k$   converges towards the varifold ${\mathbf v}_\infty$  associated to a  smooth, possibly branched, conformal minimal immersion $\vec{\Psi}_\infty$ of a constant Gauss curvature surface  $(S_\infty,h)$ equipped with a locally constant odd multiplicity $N_\infty\in C^\infty(S_\infty,2\,{\N}+1)$ and such that 
\[
\mbox{genus}\,(S_\infty)\le g\quad\mbox{ and }\quad\lim_{k\rightarrow+\infty}\mbox{Area}(\vec{\Phi}_k)=\frac{1}{2}\int_{S_\infty}N_\infty\ |d\vec{\Psi}_\infty|_h^2\ dvol_h
\]
\hfill $\Box$
\end{Th}
\medskip

 The question of comparing the {\it Morse index} of the limiting surface for the area with the {\it Morse index} of the sequence $\vec{\Phi}_k$ for the relaxed functionals $A^{\sigma_k}$ was left open in these works. The second main result of the present work is the following {\it lower semi-continuity} of the index.

\medskip
\begin{Th}
\label{th-A.3}
Let $\vec{\Phi}_k$ be a sequence of immersions of a closed surface $\Sigma$, critical points of $A^{\sigma_k}$ and such that
\[
\limsup_{k\rightarrow+\infty}A^{\sigma_k}(\vec{\Phi}_k)<+\infty\quad\mbox{and }\quad\sigma_k^2\int_{\Sigma^g}(1+|\vec{\mathbb I}_{\vec{\Phi}_k}|^2)^2\ dvol_{g_{\vec{\Phi}_k}}=o\lf( \frac{1}{\log\sigma_k^{-1}} \rg)
\]
Then there exists a subsequence $\vec{\Phi}_{k'}$ such that the corresponding immersed surface converges in varifolds towards a parametrized integer rectifiable stationary varifold ${\mathbf v}_\infty:=(S_\infty,\vec{\Psi}_\infty,N_\infty)$. If $N_\infty\equiv 1$ then we have
\be
\label{A.I.14}
\mbox{Ind } (\vec{\Psi}_\infty)\le \liminf_{k\rightarrow \infty} \mbox{Ind}^{\,\sigma_{k'}}(\vec{\Phi}_{k'})
\ee
where $\mbox{Ind} (\vec{\Psi}_\infty)$  is the maximal dimension of a subspace of $T_{[\vec{\Psi}_\infty]}{\mathfrak M}$ on which $D^2\mbox{Area}(\vec{\Psi}_\infty)$ is strictly negative and $\mbox{Ind}^{\,\sigma_k}(\vec{\Phi}_k)$ is the maximal dimension of a subspace of $T_{[\vec{\Phi}_k]}{\mathfrak M}$ on which $D^2 A^{\sigma_k}(\vec{\Phi}_k)$ is strictly negative.
\hfill $\Box$
\end{Th}
\begin{Rm}
After the present work has been completed, the author in collaboration with Alessandro Pigati proved that the condition $N_\infty\equiv 1$ always holds for the varifold
limit of sequences of critical points of $F_{\sigma_k}$ satisfying the entropy condition (\ref{entropy}) (see \cite{PiR2}). Hence the lower-semi-continuity of the index always holds.
Combining this result with the main theorem of \cite{Mic} the authors in \cite{PiR2} establish that the Morse Index of any  minimal surface realizing the minmax of a k-dimensional
homological (or co-homological) family obtained by the viscosity method is bounded by $k$. \hfill $\Box$
\end{Rm}

The paper is organized as follows. In the next section we are proving theorem~\ref{th-hilbert} from which we deduce theorem~\ref{th-1}. In a short intermediate section we establish proposition~\ref{th-fredholm}.
Then, in section 4, we are proving the lower semi continuity of the index in the {\it viscosity method} (i.e. theorem~\ref{th-A.3}).

\section{A Proof of Theorem~\ref{th-hilbert}.}
Let $\Sigma^g$ be a closed connected oriented surface of genus $g$. Let $\mbox{Diff}_+(\Sigma^g)$ be the topological group of positive  $W^{3,2}-$diffeomorphisms of $\Sigma$,  isotopic to the identity. This can be seen as an open subspace of $W^{3,2}(\Sigma,\Sigma)$ which itself defines a Hilbert Manifold (see \cite{Riv-minimax}). For $g=0$ we are marking 3 distinct points, that we denote $a_1,a_2,a_3$, for $g=1$ we are marking one point  that we denote $a$ and for $g>1$ no point is marked. We denote by $\mbox{Diff}^\ast_+(\Sigma^g)$ the sub-group of $\mbox{Diff}_+(\Sigma^g)$ which are fixing the marked points. In particular for $g>1$ we have $\mbox{Diff}^\ast_+(\Sigma^g)=\mbox{Diff}_+(\Sigma^g)$. We have the following lemma
\begin{Lm}
\label{lm-free} The action of $\mbox{Diff}^\ast_+(\Sigma^g)$ on $\mbox{Imm}(\Sigma^g,M^m)$ is free.\hfill $\Box$
\end{Lm}
\noindent{\bf Proof of Lemma~\ref{lm-free}.} We first claim that every element in $\mbox{Diff}^\ast_+(\Sigma^g)$ posses at least one fixed point. This is included in the definition for
$g=0,1$. For $g>1$ we have that for any diffeomorphism $\Psi$ isotopic to the identity the {\it Lefschetz Number} $L(\Psi)$ is given by definition by
\[
L(\Psi)=\mbox{Tr}(\Psi| H_0(\Sigma^g))-\mbox{Tr}(\Psi| H_1(\Sigma^g))+\mbox{Tr}(\Psi| H_2(\Sigma^g))=2-2g
\]
Hence for $g>1$ we have $L(\Psi)\ne 0$ and then $\Psi$ must have at least one fixed point. Due to lemma 1.3 in \cite{CMM} we deduce that for any $g\in {\N}$ the action of $\mbox{Diff}^{\, \ast}_+(\Sigma^g)$ on ${\mathcal M}^g:=\mbox{Imm}(\Sigma^g,M^m)$ is {\it free}. \hfill $\Box$

\medskip
\noindent{\bf Proof of Theorem~\ref{th-hilbert}.}
 Let $\vec{\Phi}\in \mbox{Imm}(\Sigma, M^m)$. A basis of neighborhoods  of $\vec{\Phi}$ is given by
\[
{\mathcal V}^\ep_{\vec{\Phi}}:=\lf\{\vec{\Xi}=\pi_{M^m}\lf(\vec{\Phi}+\vec{v}\rg)\quad;\quad \vec{v}\in \Gamma^{3,2}(\vec{\Phi}^\ast TM^m)\cap W^{3,2}(\Sigma,{\R}^Q)\ \mbox{ and }\ \|\vec{v}\|_{W^{3,2}}<\ep\rg\}\quad.
\]
for $\ep>0$ small enough and where $\Gamma^{3,2}(\vec{\Phi}^\ast TM^m)$ denotes the $W^{3,2}-$sections of the pull-back bundle $\vec{\Phi}^\ast TM^m$, that is the sub-vector space
of $\vec{v}\in W^{3,2}(\Sigma^g,{\R}^Q)$ such that $\vec{v}(x)\in T_{\vec{\Phi}(x)}\Sigma^g$ for any $x\in \Sigma^g$.

For any  $\vec{v}\in \Gamma^{3,2}(\vec{\Phi}^\ast TM^m)$ we consider the tensor in $\Gamma^{2,2}((T^\ast\Sigma^g)^{(0,1)}\otimes(T\Sigma^g)^{(1,0)})$ given by
\[
  \ov{D}^\ast_{\vec{\Phi}}\vec{v}\ \res\, g^{-1}_{\vec{\Phi}} \quad\mbox{ where }\quad g_{\vec{\Phi}}^{-1}=e^{-2\la}\ [\p_z\otimes\p_{\ov{z}}+\p_{\ov{z}}\otimes \p_z]
\]
where
\[
 \ov{D}^\ast_{\vec{\Phi}}\vec{v}=\ov{P}_{\vec{\Phi}}^\perp\lf(\ov{\p}\vec{w}\ \dot{\otimes}\ \ov{\p}\vec{\Phi}\rg)
\]
$\ov{P}_{\vec{\Phi}}$ the $L^2$ projection orthogonal projection from $(\wedge^{(0,1)}\Sigma)^{\otimes^2}$ onto the space  of {\it anti-holomorphic quadratic forms} $\mbox{AHol}_Q(\Sigma,g_{\vec{\Phi}})$ on $(\Sigma,g_{\vec{\Phi}})$  and by $\ov{P}_{\vec{\Phi}}^\perp:=\mbox{Id}-\ov{P}_{\vec{\Phi}}$.
and where $\res$ is the contraction operator between covariant and contravariant tensors. In particular we have in local complex coordinates
\[
\lf(b\ d{\ov{ z}}\otimes d{\ov{ z}}\rg)\res g^{-1}_{\vec{\Phi}}  =e^{-2\la}\ b\ d{\ov{ z}}\otimes \p_z
\]
We denote
\be
\label{mathcali}
{\mathcal I}:=\lf\{ \ov{D}^\ast_{\vec{\Phi}}\vec{v}\ \res\, g^{-1}_{\vec{\Phi}} \ ;\quad \vec{v}\in \Gamma^{3,2}(\vec{\Phi}^\ast TM^m)\rg\}\quad.
\ee
Recall the definition of the $\ov{\p}$ operator on $\wedge^{(1,0)}\Sigma^g$ given in local coordinates by
\[
\ov{\p}\lf(a\,\p_z\rg)=\p_{\ov{z}}a\ d\ov{z}\otimes\p_z\quad.
\]
Denote $\mbox{Hol}_1(\Sigma^g)$ the finite dimensional subspace of $\Gamma^{3,2}(\wedge^{(1,0)}\Sigma^g)$ made of holomorphic sections\footnote{Due to Riemann-Hurwitz Theorem, the  holomorphic tangent bundle $T^{(1,0)}\Sigma^g$ which is the inverse of the canonical bundle of the Riemann surface defined by $(\Sigma,g_{\vec{\Phi}})$ has a degree given by
\[
\mbox{deg}\lf(T^{(1,0)}\Sigma^g\rg)=\int_{\Sigma^g}c_1\lf(T^{(1,0)}\Sigma^g\rg)=2-2g
\]
 and therefore $\mbox{Hol}_1(\Sigma^g)\ne 0$ $\Rightarrow$ $g<2$.} of $T^{(1,0)}\Sigma^g$. 
We shall now prove the following lemma.

\begin{Lm}
\label{lm-delbar}
Under the previous notations we have that $\forall\ \vec{v}\in W^{3,2}(\Sigma^g,{\R}^Q)$
\be
\label{delbar}
\forall\ \vec{v}\in W^{3,2}(\Sigma,{\R}^Q)\quad \exists\, !\, f\in (\mbox{Hol}_1(\Sigma^g))^\perp\cap \Gamma^{3,2}(\wedge^{(1,0)}\Sigma^g)\quad\quad\mbox{s. t. }\quad\ov{\p}f=\ov{D}^\ast_{\vec{\Phi}}\vec{v}\ \res\, g^{-1}_{\vec{\Phi}}\quad.
\ee
Moreover
\be
\label{delbar-estim}
\|f\|_{W^{3,2}}\le C_{\vec{\Phi}}\ \|\vec{v}\|_{W^{3,2}}\quad.
\ee
\hfill $\Box$
\end{Lm}
\noindent{\bf Proof of lemma~\ref{lm-delbar}.}
We have for any $\al=a\, \p_z\in \Gamma^{3,2}(\wedge^{(1,0)}\Sigma^g)$ and $\beta=b\ d\ov{z}\otimes\p_z\in \Gamma^{2,2}((T^\ast\Sigma^g)^{(0,1)}\otimes(T\Sigma^g)^{(1,0)})$
\[
\begin{array}{l}
\ds\int_{\Sigma^g}\lf<\ov{\p}\lf(a\,\p_z\rg), b\ d\ov{z}\otimes \p_z\rg>_{g_{\vec{\Phi}}}\ d\mbox{vol}_{g_{\vec{\Phi}}}=\Re\lf[\frac{i}{2}\int_{\Sigma^g} {\p}\ov{a}\  b\  e^{2\la}\, dz\wedge d\ov{z}\rg]\\[5mm]
\ds\quad=\Re\lf[\frac{i}{2}\int_{\Sigma^g} d[\ov{a}\  b\  e^{2\la}]\wedge d\ov{z}\rg]-\Re\lf[\frac{i}{2}\int_{\Sigma^g} \ov{a}\  {\p}_z[b\  e^{2\la}]\ dz\wedge d\ov{z}\rg]\\[5mm]
\ds\quad=\int_{\Sigma^g} \lf<\al,\lf(\p\lf(\beta\res g_{\vec{\Phi}}\rg)\res_2\, g_{\vec{\Phi}}^{-1}\rg)\res g_{\vec{\Phi}}^{-1}\rg>_{g_{\vec{\Phi}}}\ d\mbox{vol}_{g_{\vec{\Phi}}}
\end{array}
\]
where
\[
g_{\vec{\Phi}}^{-1}=e^{-2\la}\ [\p_z\otimes\p_{\ov{z}}+\p_{\ov{z}}\otimes \p_z]\quad\mbox{ , }\quad (b \ dz\otimes d\ov{z}\otimes d\ov{z})\res_2 g_{\vec{\Phi}}^{-1}=e^{-2\la}\, b\ d\ov{z}\quad\mbox{and }\quad (e^{-2\la}\, b\ d\ov{z})\res g_{\vec{\Phi}}^{-1}=e^{-4\la}\, b\ \p_z
\]
Hence we have proved that the adjoint of $\ov{\p}$ on $\Gamma((T^\ast\Sigma^g)^{(0,1)}\otimes(T\Sigma^g)^{(1,0)})$ is given by
\[
\ov{\p}^\ast\ :\ \beta\in \Gamma((T^\ast\Sigma^g)^{(0,1)}\otimes(T\Sigma^g)^{(1,0)})\ \longrightarrow\ \ov{\p}^\ast\beta=\lf(\p\lf(\beta\res g_{\vec{\Phi}}\rg)\res_2\, g_{\vec{\Phi}}^{-1}\rg)\res g_{\vec{\Phi}}^{-1}\in \Gamma(\wedge^{(1,0)}\Sigma^g)
\]
We have $\mbox{Im}\,\ov{\p}=(\mbox{Ker}\,\ov{\p}^\ast)^\perp$. We have that 
\be
\label{kerdelbarstar}
\mbox{Ker}\,\ov{\p}^\ast=\lf\{\beta\in \Gamma((T^\ast\Sigma^g)^{(0,1)}\otimes(T\Sigma^g)^{(1,0)})\quad;\quad\beta\res g_{\vec{\Phi}}\in \mbox{Hol}_Q(\Sigma^g,g_{\vec{\Phi}}))\rg\}
\ee
Observe that
\[
\forall\,\gamma\in \Gamma\lf(\lf((T^\ast\Sigma^g)^{(0,1)}\rg)^{\otimes^2}\rg)\ ,\ \forall\,\beta\in \Gamma((T^\ast\Sigma^g)^{(0,1)}\otimes(T\Sigma^g)^{(1,0)})\quad
\lf<\gamma\res g^{-1}_{\vec{\Phi}},\beta\rg>_{g_{\vec{\Phi}}}=\lf<\gamma,\beta\res g_{\vec{\Phi}}\rg>_{g_{\vec{\Phi}}}
\]
This implies 
\be
\label{carac}
\gamma\res g^{-1}_{\vec{\Phi}}\in \mbox{Im}\,\ov{\p}\quad\Longleftrightarrow\quad \gamma\in \lf(\mbox{AHol}_Q(\Sigma^g,g_{\vec{\Phi}})\rg)^\perp\quad.
\ee
We deduce (\ref{delbar}) from (\ref{carac}) and (\ref{delbar-estim}) follows by classical estimates for elliptic complexes in Sobolev Spaces.\hfill $\Box$

\medskip

\noindent{\bf Continuation of the proof of theorem~\ref{th-hilbert}.} To $f\in(\mbox{Hol}_1(\Sigma^g))^\perp\cap \Gamma^{3,2}(\wedge^{(1,0)}\Sigma^g)$ solving (\ref{delbar}) we assign
\[
X:=2\,\Re(f)=2\,\Re((f_1+i\,f_2)\ \p_z)=(f_1\,\p_{x_1}+f_2\,\p_{x_2})=X_1\, \p_{x_1}+X_2\,\p_{x_2}\quad.
\]
Observe that, if we denote $\vec{X}:=d\vec{\Phi}\cdot X$, we have
\[
\ov{\p}\lf(\vec{X}\cdot \ov{\p}\vec{\Phi}\res g_{\vec{\Phi}}^{-1}\rg)=\ov{\p}\lf(e^{2\la} (X_1+i\,X_2) \ d\ov{z}\res g_{\vec{\Phi}}^{-1}\rg)=\ov{\p} f
\]
Observe also that, since $\vec{X}$ is tangent to the immersion $\vec{X}\cdot\ov{\p}\lf( \ov{\p}\vec{\Phi}\res g_{\vec{\Phi}}^{-1}\rg)=0$ hence
\[
\ov{\p}f=(\ov{\p}\vec{X}\cdot\ov{\p}\vec{\Phi})\res g_{\vec{\Phi}}^{-1}
\]
Using $\mbox{Im}\,\ov{\p}=(\mbox{Ker}\,\ov{\p}^\ast)^\perp$ and the characterization of $\mbox{Ker}\,\ov{\p}^\ast$ given by (\ref{kerdelbarstar}), we have $$\ov{\p}\vec{X}\cdot\ov{\p}\vec{\Phi}=P^\perp_{\vec{\Phi}}(\ov{\p}\vec{X}\cdot\ov{\p}\vec{\Phi})$$ and hence
\be
\label{rep}
\ov{\p}f=\ov{D}_{\vec{\Phi}}^\ast\vec{X} \res g_{\vec{\Phi}}^{-1}
\ee
We denote
\[
\lf\{
\begin{array}{l}
\ds{\mathcal X}^{3,2}(S^2):=\lf\{{X}\in \Gamma^{3,2}(TS^2)\ ;\ X(a_i)=0\quad i=1,2,3\rg\}\quad,\\[3mm]
\ds{\mathcal X}^{3,2}(T^2):=\lf\{{X}\in \Gamma^{3,2}(T T^2)\ ;\ X(a)=0\rg\}\quad,\\[3mm]
\ds{\mathcal X}^{3,2}(\Sigma^g)=\Gamma^{3,2}(\Sigma^g)\quad\mbox{ for }g>1
\end{array}
\rg.
\]
The space of Holomorphic Vector Field on $T^{(1,0)}S^2$ is a $3-$dimensional complex vector space given in ${\C}$, after the stereographic projection, by
\[
h(z)= (\al+\beta\, z+\gamma\, z^2)\ \p_z\quad\mbox{ where }(\al,\beta,\gamma)\in {\C}^3\quad.
\]
Whereas the space of Holomorphic Vector Field on $T^{(1,0)}T^2$ is a $1-$dimensional complex vector space given in ${\C}$ by
\[
h(z)=\al\ \p_z\quad\mbox{ where }\al\in {\C}\quad.
\]
while for $g>1$ we have $\mbox{Hol}_1(\Sigma^g)=\{0\}$. Hence, for any $g\in {\N}$ and any
\be
\label{hol-1}
f\in(\mbox{Hol}_1(\Sigma^g))^\perp\cap \Gamma^{3,2}(\wedge^{(1,0)}\Sigma^g)\quad\exists\ !\
h_f\in \mbox{Hol}_1(\Sigma^g)\quad\mbox{ s. t. }\quad\Re(f+h_f)\in {\mathcal X}^{3,2}(\Sigma^g)
\ee
moreover the map $f\rightarrow h_f$ from $(\mbox{Hol}_1(\Sigma^g))^\perp\cap \Gamma^{3,2}(\wedge^{(1,0)}\Sigma^g)$ into $\mbox{Hol}_1(\Sigma^g)$ is linear and smooth.
Hence we can summarize what we have proved so far in the following lemma.
\begin{Lm}
\label{lm-decomposition} Let $\vec{\Phi}$ be a $W^{3,2}-$immersion. Then the following holds
\[
\begin{array}{l}
\ds\forall\, \vec{v}\in \Gamma^{3,2}(\vec{\Phi}^\ast TM^m)\quad\exists \ !\ X\in{\mathcal X}^{3,2}(\Sigma^g) \quad\mbox{ s. t. }\quad\\[3mm]
\ds\ov{\p}\lf(X-\,i\,X^\perp\rg)=\ov{D}_{\vec{\Phi}}^\ast\vec{X}\res g_{\vec{\Phi}}^{-1}=\ov{D}_{\vec{\Phi}}^\ast\vec{v}\res g_{\vec{\Phi}}^{-1}
\end{array}
\]
where $\vec{X}=d\vec{\Phi}\cdot X$ and such that
\[
\|{X}\|_{W^{3,2}}\le C_{\vec{\Phi}}\ \|\vec{v}\|_{W^{3,2}}
\]
\hfill $\Box$
\end{Lm}
\noindent{\bf End of the proof of theorem~\ref{th-hilbert}.} Let $g_0$ be a smooth reference metric on $\Sigma^g$ and denote by $\exp^{g_0}$ the smooth
exponential map from $T\Sigma$ into $\Sigma$ associated to $g_0$. Let $\ep>0$ small and denote respectively
\[
{\mathcal X}^{3,2}_\ep(\Sigma^g):=\lf\{{X}\in {\mathcal X}^{3,2}(\Sigma^g)\ ;\ \|X\|_{W^{3,2}}<\ep\rg\}\quad,
\]
and
\[
{\mathcal D}^\perp_\ep:=\lf\{\Psi\in \mbox{Diff}_+^\ast(\Sigma)\ ;\quad \exists\, X\in {\mathcal X}^{3,2}_\ep(\Sigma^g)\quad \mbox{s.t}\quad \Psi(x)= \exp^{g_0}_x(X(x)) \rg\}\quad.
\]
We define
\[
\begin{array}{rcl}
\Lambda_{\vec{\Phi}}\ :\ {\mathcal V}_{\vec{\Phi}}^\ep\times {\mathcal D}^\perp_\ep &\longrightarrow& \Gamma^{2,2}((T^\ast\Sigma )^{(0,1)}\otimes (T^\ast\Sigma )^{(0,1)})\\[3mm]
(\vec{\Xi},\Psi) &\longrightarrow& \ov{D}_{\vec{\Phi}}^\ast\lf(\vec{\Xi}\circ\Psi\rg)\res g^{-1}_{\vec{\Phi}}
\end{array}
\]
The map is clearly $C^1$ and lemma~\ref{lm-decomposition}  gives that
\[
\lf.\p_{\Psi}\Lambda_{\vec{\Phi}}\rg|_{(\vec{\Phi},0)}\cdot X=\ov{D}_{\vec{\Phi}}^\ast\lf(d\vec{\Phi}\cdot X\rg)\res g^{-1}_{\vec{\Phi}}
\]
realizes an isomorphism between ${\mathcal X}^{3,2}$ and ${\mathcal I}$ (defined in (\ref{mathcali})). The implicit function theorem gives then the existence
of a $C^1$ map $\Psi_{\vec{\Phi}}(\vec{\Xi})$ defined in a neighborhood of $\vec{\Phi}$ such that
\[
\ov{D}_{\vec{\Phi}}^\ast\lf(\vec{\Xi}\circ\Psi_{\vec{\Phi}}(\vec{\Xi})\rg)\res g^{-1}_{\vec{\Phi}}=0
\]
and we denote $\vec{w}_{\vec{\Phi}}(\vec{\Xi}):=\vec{\Xi}\circ\Psi_{\vec{\Phi}}(\vec{\Xi})-\vec{\Phi}$.

\medskip

For any element $\Psi_0\in \mbox{Diff}_+^\ast(\Sigma)$ close to the identity and $\vec{\Xi}$ close enough to $\vec{\Phi}$ one has trivially
\[
\ov{D}_{\vec{\Phi}}^\ast\lf(\vec{\Xi}\circ\Psi_0\circ\Psi_0^{-1}\circ\Psi_{\vec{\Phi}}(\vec{\Xi})\rg)\res g^{-1}_{\vec{\Phi}}=0\quad.
\]

 Because of the local uniqueness of $\Psi_{\vec{\Phi}}(\vec{\Xi})$ given by
the implicit function theorem, we deduce the equivariance property
\[
\Psi_{\vec{\Phi}}(\vec{\Xi}\circ{\Psi}_0)=\Psi_0^{-1}\circ\Psi_{\vec{\Phi}}(\vec{\Xi})\quad\mbox{ and }\quad\vec{w}_{\vec{\Phi}}(\vec{\Xi}\circ\Psi_0):=\vec{\Xi}\circ\Psi_{\vec{\Phi}}(\vec{\Xi})-\vec{\Phi}=\vec{w}_{\vec{\Phi}}(\vec{\Xi})
\]
We then naturally extend, by equivariance, the map ${\mathcal T}_{\vec{\Phi}}:=(\vec{w}_{\vec{\Phi}},\Psi_{\vec{\Phi}})$ on a neighborhood ${\mathcal O}_{\vec{\Phi}}$ of $\vec{\Phi}$ invariant under the action of $\mbox{Diff}_+^\ast(\Sigma)$.

\medskip

We are now proving the Hausdorff property for ${\mathfrak M}_{g}(\Sigma^g,M^m):=\mbox{Imm}_0(\Sigma^g,M^m)/\mbox{Diff}^{\, \ast}_+(\Sigma^g)$. Following classical considerations
(see the arguments in the proof of lemma 2.9.9 of \cite{Var}) it suffices to prove that 
\[
\Gamma:=\lf\{ (\vec{\Phi},\vec{\Phi}\circ\Psi)\quad;\quad\vec{\Phi}\in \mbox{Imm}_0(\Sigma^g,M^m)\ \mbox{ and }\ \Psi\in \mbox{Diff}^{\, \ast}_+(\Sigma^g)\rg\}
\]
is closed in $(\mbox{Imm}_0(\Sigma^g,M^m))^2$. This follows fro the first part of the proof of the theorem. Let $(\vec{\Phi}_k,\vec{\Xi}_k:=\vec{\Phi}_k\circ\Psi_k)\rightarrow (\vec{\Phi}_\infty,\vec{\Xi}_\infty)$ in $W^{3,2}$. For $k$ large enough both $\vec{\Phi}_k$ and $\vec{\Xi}_k$ are included in ${\mathcal O}_{\vec{\Phi}_\infty}$. Because of the continuity of the map
 $\vec{w}_{\vec{\Phi}_\infty}$ we have respectively
 \[
 \vec{w}_{\vec{\Phi}_\infty}(\vec{\Phi}_k)\rightarrow  \vec{w}_{\vec{\Phi}_\infty}(\vec{\Phi}_\infty)=0\quad\mbox{ and }\quad\vec{w}_{\vec{\Phi}_\infty}(\vec{\Xi}_k)\rightarrow  \vec{w}_{\vec{\Phi}_\infty}(\vec{\Xi}_\infty)
 \]
The equivariance of $\vec{w}_{\vec{\Phi}_\infty}$ gives $\vec{w}_{\vec{\Phi}_\infty}(\vec{\Xi}_k)=\vec{w}_{\vec{\Phi}_\infty}(\vec{\Phi}_k)$ hence  $\vec{w}_{\vec{\Phi}_\infty}(\vec{\Xi}_\infty)=0$.
Thus $\vec{\Xi}_\infty\circ\Psi_{\vec{\Phi}_\infty}=\vec{\Phi}_\infty$ and this shows that $\Gamma$ is closed and then ${\mathfrak M}(\Sigma^g,M^m)$ defines an Hausdorff {\it Hilbert Manifold}
and theorem~\ref{th-hilbert} is proved. \hfill $\Box$
\section{A Proof of Proposition~\ref{th-fredholm}.}
From \cite{KLL} (see also an alternative approach in \cite{BR}) we know that under the assumptions that $\vec{\Phi}$ is a critical point of $F$, it defines a smooth immersion in conformal coordinates. We shall be working in the chart in the neighborhood
of $[\vec{\Phi}]$ in ${\mathfrak M}(\Sigma^g,M^m)$ given
by $\vec{w}_{\vec{\Phi}}$ from theorem~\ref{th-hilbert}. In other words  we identify 

\be
\label{fred-1}
T_{[\vec{\Phi}]}{\mathfrak M}\simeq\lf\{\vec{w}\in \Gamma^{3,2}(\vec{\Phi}^\ast T M^m)\quad;\quad P^\perp_{\vec{\Phi}}\lf(\ov{\p}\vec{w}\,\dot{\otimes}\,\ov{\p}\vec{\Phi}\rg)=0\rg\}
\ee
For such a $\vec{w}$ we denote by $q_{\vec{w}}$ the holomorphic quadratic form given by
\[
\ov{\p}\vec{w}\,\dot{\otimes}\,\ov{\p}\vec{\Phi}=q_{\vec{w}}\quad.
\]
After contracting with the tensor $g_{\vec{\Phi}}^{-1}$, this equation becomes
\[
\ov{\p}\lf(\vec{w}\cdot\ov{\p}\vec{\Phi}\res g_{\vec{\Phi}}^{-1}\rg)=-\pi_{\vec{n}}(\vec{w})\cdot \vec{h}^0+q_{\vec{w}}\res g_{\vec{\Phi}}^{-1}
\]
where $\vec{h}^0$ is the trace free part of the second fundamental form, which is orthogonal to the tangent plane of the immersion and given in local coordinates by 
\[
\vec{h}^0_{\vec{\Phi}}=\p_{\ov{z}}\lf(e^{-2\la}\p_{\ov{z}}\vec{\Phi}\rg)\ d\ov{z}\otimes \p_z\quad.
\]
Using the characterization of $\mbox{Im}\,\ov{\p}=(\mbox{Ker}\,\ov{\p}^\ast)^\perp$ given by  (\ref{carac}) we deduce
\[
\ov{\p}\lf(\vec{w}\cdot\ov{\p}\vec{\Phi}\res g_{\vec{\Phi}}^{-1}\rg)=-{\mathfrak P}_{\vec{\Phi}}\lf(\pi_{\vec{n}}(\vec{w})\cdot \vec{h}^0_{\vec{\Phi}}\rg)
\]
where ${\mathfrak P}_{\vec{\Phi}}$ is the orthogonal projection onto $(\mbox{Hol}_Q(\Sigma^g,g_{\vec{\Phi}})\res  g_{\vec{\Phi}}^{-1})^\perp$. Denote $\vec{X}_{\vec{w}}$ the projection
of $\vec{w}$ onto the tangent plane (i.e. $\vec{X}_{\vec{w}}=\vec{w}-\pi_{\vec{n}}(\vec{w})$) and let $X_{\vec{w}}$ be the vector field on $\Sigma$ such that $d\vec{\Phi}\cdot X_{\vec{w}}=\vec{X}_{\vec{w}}$. Following the computations from the previous subsection we deduce
\be
\label{fred-2}
\ov{\p}\lf(X_{\vec{w}}-\,i\,X_{\vec{w}}^\perp\rg)=-{\mathfrak P}_{\vec{\Phi}}\lf(\pi_{\vec{n}}(\vec{w})\cdot \vec{h}^0_{\vec{\Phi}}\rg)
\ee
Denote
\[
\begin{array}{rcl}
\ds\pi_T\ :\ \Gamma^{3,2}(\vec{\Phi}^\ast TM^m) &\longrightarrow & \Gamma^{3,2}((T\Sigma)^{(1,0)})\\[3mm]
\ds\vec{w} &\longrightarrow & X_{\vec{w}}-\,i\,X_{\vec{w}}^\perp
\end{array}
\]
In view of the expression of the second derivative $D^2F$ given by (\ref{deuxi-18}) we have that, modulo compact operators (remembering that $\vec{\Phi}$ is smooth), We are reduced\footnote{The sum of a {\it Fredholm operator} with a {\it compact operator} is {\it Fredholm}.} to study the Fredholm nature of the operator generated by the following quadratic form
\[
\begin{array}{l}
\ds Q_{\vec{\Phi}}(\vec{w})= \,\, \int_{\Sigma}(1+|\vec{\mathbb I}_{\vec{\Phi}}|^2_{g_{\vec{\Phi}}}) \lf|\pi_{\vec{n}}\lf(D^{g_{\vec{\Phi}}}d\vec{w}\rg)\rg|_{g_{\vec{\Phi}}}^2\ dvol_{g_{\vec{\Phi}}}
+2\, \int_{\Sigma}\, \lf|\lf<\vec{\mathbb I}_{\vec{\Phi}},D^{g_{\vec{\Phi}}}d\vec{w}\rg>_{g_{\vec{\Phi}}}\rg|^2\ dvol_{g_{\vec{\Phi}}}
\end{array}
\]
combined with (\ref{fred-2}). Hence the Symbols of the generated operator, in local conformal coordinates, is given by
\[
\lf\{
\begin{array}{l}
2\,e^{-2\la} \pi_{\vec{n}}\circ\lf[(1+|\vec{\mathbb I}_{\vec{\Phi}}|^2_{g_{\vec{\Phi}}}) |\xi|^4+2\,e^{-4\la}\sum_{i,j,k,l} \vec{\mathbb I}_{kl}\otimes\vec{\mathbb I}_{kl} \ \xi_i\,\xi_j\,\xi_k\,\xi_l\rg]\circ\pi_{\vec{n}}\\[5mm]
(\xi_1+i\,\xi_2)\circ\pi_T
\end{array}
\rg.
\]
This is clearly the symbol defining an {\it elliptic operator} on $\Gamma^{3,2}(\vec{\Phi}^\ast TM^m)$ and $D^2F$ is {\it Fredholm} on $T_{[\vec{\Phi}]}{\mathfrak M}$. This concludes the proof
of proposition~\ref{th-fredholm}.
\section{A Proof of Theorem~\ref{th-A.3}, the Lower Semi Continuity of the Index.}
We shall assume that $\Sigma^g$ is connected. We shall present the computations  for $M^m=S^m$. The general constraint generates lower order terms whose abundance could mask the true reason why the theorem is true whereas the same terms in the $M^m=S^m$ case are easier to present. The first part of the theorem is the  main results of \cite{Riv-minmax} . It remains to prove the inequality
(\ref{A.I.14}) under the assumption of theorem~\ref{th-A.3}. The first derivative of the area of an immersion (possibly branched) of a closed surface $\Sigma$ into ${\R}^Q$ is given by (see \cite{Riv-minmax})
\[
D\mbox{Area}(\vec{\Phi})\cdot\vec{w}=\int_\Sigma \lf<d\vec{\Phi}\,;\,d\vec{w}\rg>_{g_{\vec{\Phi}}}\ d\mbox{vol}_{g_{\vec{\Phi}}}
\]
and the second derivative\footnote{ A reader familiar to the rich literature in geometry  on minimal surface theory in 3 dimension might not immediately recognize the most commonly used expression of the second derivative of the area by the mean of the {\it Jacobi field}. This classical presentation of $D^2\mbox{Area}$ has the advantage to ``reduce'' this operator to an operator on function by introducing the decomposition $\vec{w}=w\, \vec{n}$. This decomposition however is not ``analytically'' favorable since $\vec{n}$ has a-priori one degree of regularity less than $\vec{w}$. This observation is at the base of the analysis of the {\it Willmore functional} as it has been developed by the author in the recent years.   }
\[
D^2 \mbox{Area}(\vec{\Phi})\cdot(\vec{w},\vec{w})=\int_\Sigma \lf[\lf<d\vec{w}\,;\,d\vec{w}\rg>_{g_{\vec{\Phi}}}+\lf|\lf<d\vec{\Phi}\,;\,d\vec{w}\rg>_{g_{\vec{\Phi}}}\rg|^2 - 2^{-1}\lf|d\vec{\Phi}\dot{\otimes}\, d\vec{w}+d\vec{w}\dot{\otimes}\, d\vec{\Phi}\rg|^2\rg]\ d\mbox{vol}_{g_{\vec{\Phi}}}
\]
where in coordinates $d\vec{\Phi}\dot{\otimes} d\vec{w}:=\sum_{i,j}\p_{x_i}\vec{\Phi}\cdot\p_{x_j}\vec{w}\ dx_i\otimes dx_j$. 

\medskip

Since we are assuming $N_\infty\equiv1$  a.e. on $S_\infty$ and, following the proof of the main theorem of \cite{Riv-minmax}, 
we can extract a subsequence that we keep denoting $\vec{\Phi}_{k}$ such that we have a {\it bubble tree } strong $W^{1,2}$ convergence
of $\vec{\Phi}_{k}$ towards a minimal (possibly branched) immersion $\vec{\Psi}_\infty$ of a surface $S_\infty$. More precisely, if one denotes $\{S^j_\infty\}_{j\in J}$ to be the connected components of $S_\infty$, for every $j\in J$ there exists $N^j$ points $\{a^{j,l}\}_{l=1\cdots N^j}$ (containing in particular the possible branched points of $\vec{\Psi}_\infty$ and a converging sequence of constant scalar curvature metrics  $h^j_{k}$ of volume one  and for any $\delta>0$ a sequence of conformal embeddings $\phi^j_{k}$ from $(S^j_\infty\setminus \cup_{l=1}^{N^j} B_\delta(a^{j,l}), h^j_{k})$ into $(\Sigma^g,g_{\vec{\Phi}_{k}})$ such that
\be
\label{A.I.15}
\vec{\Psi}^j_{k}:=\vec{\Phi}_{k}\circ \phi^j_{k}\longrightarrow {\vec{\Psi}_\infty}\quad\quad\mbox{ strongly in } W^{1,2}_{loc}(S^j_\infty\setminus \cup_{l=1}^{N^j} B_\delta(a^{j,l}))
\ee

For $\delta$ small enough and $k'$ large enough the  subdomains $\Omega_{k}^j(\delta):=\phi^j_{k}(S^j_\infty\setminus \cup_{l=1}^{N^j} B_\delta(a^{j,l}))$ are disjoint and 
\[
\lim_{\delta\rightarrow 0}\lim_{k\rightarrow +\infty}\mbox{Area}\lf(\vec{\Phi}_{k}\lf(\Sigma^g\setminus \bigcup_{j\in J}\Omega_{k}^j(\delta)\rg)\rg)=0
\]
Let $\vec{w}_1\cdots\vec{w}_N$ a family of $N$ independent smooth vectors in $W^{3,2}(\vec{\Psi}_\infty^\ast T{M}^m)$ representing $N$ independent directions in $T_{[\vec{\Psi}_\infty]}{\mathfrak M}$ on the Span of which
$D^2\mbox{Area}$ is strictly negative. We can assume without loss of generality that the $\vec{w}_i$ are $C^\infty$. One modifies each of these vectors in the following way. For each $i\in\{1\cdots Q\}$ for each $j\in J$ and each $l\in\{1\cdots N^j\}$ we introduce (after identifying for each  $j$ and $l$ the tangent planes to $M^m$ around $\vec{\Phi}_\infty(a^{j,l})$ with the one at exactly $\vec{\Phi}_\infty(a^{j,l})$)
\[
\vec{w}^{\,\delta}_i(x)=\lf\{\begin{array}{l}
\vec{w}_i(x)\quad\quad\quad\quad \mbox{for }|a^{j,l}-x|\ge \sqrt{\delta}\\[3mm]
\vec{w}_i(x)\ \chi^{\,\delta}(|x-a^{j,l}|)\quad \quad\mbox{ for }\delta\le |a^{j,l}-x|\le \sqrt{\delta}\\[3mm]
0\quad\quad\quad\quad \mbox{ for } |a^{j,l}-x|\le \delta
\end{array}
\rg.
\] 
 where we take $\chi^{\,\delta}(s)$ to be  a slight smoothing of
$
\log({s}/{\delta})/\log({1}/{\sqrt{\delta}})
$.
A short computation gives that
 \[
\vec{w}_i^{\,\delta}\longrightarrow\vec{w}_i\quad\mbox{ strongly in }\quad W^{1,2}(S_\infty,{\R}^Q)\quad.
\] 
Therefore, in view of the explicit expression of $D^2 \mbox{Area}(\vec{\Psi}_\infty)\cdot(\vec{w},\vec{w})$, there exists $\delta$ small enough such that $\vec{w}^{\,\delta}_1\cdots\vec{w}^{\,\delta}_N$ realizes a family of $N$ independent smooth vectors in $W^{3,2}(\vec{\Psi}_\infty^\ast TM^m)$ on the Span of which
$D^2\mbox{Area}$ is strictly negative. We fix such a $\delta$.

Let $\rho >0$ small enough such that for any $z\in M^m$ the map $\vec{\Psi}_\infty$ is injective
on each components of $\vec{\Psi}^{-1}_\infty(\ov{B^Q_\rho(z)})\subset S^j_\infty$. Let $\{\chi_s(z)\}_{s\in \{1\cdots N\}}$ be a finite smooth partition of unity of $M^m\subset {\R}^Q$ such that the support of every $\chi_s$ is included in an $m-$ball of radius $\rho$. We denote the connected components of $\vec{\Psi}_\infty^{-1}(\mbox{Supp}(\chi_s))$ in $S_\infty$ by $\Om_{s}^t$ for $t=1\cdots n_{s}$ and $\om_s^t$ are the corresponding characteristic functions. We have that $\chi_s(\vec{\Psi}_\infty(x))\ \om_s^t(x)$ is smooth for any $s\in \{1\cdots N\}$ and any
 $t\in\{1\cdots n_s\}$ and moreover
 \[
 d(\chi_s(\vec{\Psi}_\infty(x))\ \om_s^t(x))=d(\chi_s(\vec{\Psi}_\infty(x)))\ \om_s^t(x)
\]
 We can then write each $\vec{w}_i^{\,\delta}$ in the form
\[
\vec{w}_i^{\,\delta}(x)=\sum_{s=1}^N\chi_s(\vec{\Psi}_\infty(x))\sum_{t=1}^{n_{s}}\vec{v}_{i,s}^{\,t}(\vec{\Psi}_\infty(x))\ \om_s^t
\]
where $\vec{v}^{\,t}_{i,s}$ are smooth functions\footnote{This is due to the fact that $\vec{\Psi}_\infty$ is smooth embedding on each open set $\Om_{s}^t$}. For any $s=\in \{1\cdots N\}$ since the components $\ov{\Om_s^t}$ are disjoint to each other for $t\in\{1\cdots n_s\}$ We can include them in strictly larger disjoint open sets $\ov{\Om_s^t}\subset\ti{\Om}_s^t$ moreover, because of the strong $W^{1,2}-$convergence of $\vec{\Phi}_k$ towards $\vec{\Phi}_\infty$ in $S^j_\infty\setminus \cup_{l=1}^{N^j} B_\delta(a^{j,l})$  
\be
\label{conv-bord}
\|\vec{\Psi}_k-\vec{\Psi}_\infty\|_{L^\infty(\p \ti{\Om_s^t})}\longrightarrow 0
\ee
We denote $\ti{\om}_s^t$ the  characteristic functions of $\ti{\Om_s^t}$: $\ti{\om}_s^t:={\mathbf 1}_{\ti{\Om_s^t}}$. We still have of course
\[
\vec{w}_i^{\,\delta}(x)=\sum_{s=1}^N\chi_s(\vec{\Psi}_\infty(x))\sum_{t=1}^{n_{s}}\vec{v}_{i,s}^{\,t}(\vec{\Psi}_\infty(x))\ \ti{\om}_s^t
\]
Because of (\ref{conv-bord}), we have that
\[
\mbox{dist}\lf(\vec{\Psi}_k(\p \ti{\Om_s^t}),\vec{\Psi}_\infty(\p \ti{\Om}_s^t\rg)\longrightarrow 0
\]
Since $\chi_s$ is zero in an open neighborhood of each $\vec{\Psi}_\infty(\p \ti{\Om}_s^t)$ and since $\vec{\Psi}_k$ is smooth, we have  for every
$s$ and $t$ and for $k$ large enough
\be
\label{canc}
\chi_s\circ\vec{\Psi}_k\equiv 0 \quad \mbox{ in an open neighborhood }U_{k,s}^t\mbox{ of } \p \ti{\Om}_s^t\quad.
\ee
Hence in particular we have for every
$s$ and $t$ and $k$ large enough
\[
 d(\chi_s(\vec{\Psi}_k(x))\ \tilde{\om}_s^t(x))=d(\chi_s(\vec{\Psi}_k(x)))\ \tilde{\om}_s^t(x)
\]
It is then clear that 
\be
\label{A.I.16-a}
\vec{w}_{i,k}^{\,\delta}(x):=\sum_{s=1}^N\chi_s(\vec{\Psi}_k(x))\sum_{t=1}^{n_{s}}\vec{v}_{t,s}(\vec{\Psi}_k(x))\ \ti{\om}_s^t\ \longrightarrow\ \vec{w}_i^{\,\delta}(x)\mbox{ strongly in }W^{1,2}_{loc}(S^j_\infty\setminus \cup_{l=1}^{N^j} B_\delta(a^{j,l}))
\ee
Using the compositions with the  maps $(\phi^{j,k'})^{-1}$ we extend the $\vec{w}^{\,\delta}_{i,k}$, that we still denote $\vec{w}_{i,k}^{\,\delta}$ to the whole of $\Sigma^g$ by taking
$\vec{w}_{i,k}^{\,\delta}=0$ on $\Sigma^g\setminus \bigcup_{j\in J}\Omega_{k'}^j(\delta)$.  We see $\vec{w}^{\,\delta}_{i,k}$ as vectors in ${\R}^Q$ and we denote by $\pi^j_{k'}$ the map from $S^j_\infty\setminus \cup_{l=1}^{N^j} B_\delta(a^{j,l})$ into the space of projection matrices which to $x\in S^j_\infty\setminus \cup_{l=1}^{N^j} B_\delta(a^{j,l})$ assigns the orthogonal projection from $T_{\vec{\Psi}^j_{k'}(x)}{\R}^Q$ into $T_{\vec{\Psi}^j_{k'}(x)}M^m$. In other words, let $P_z$ to be the $C^1$ map from $M^m$ into the space of $Q\times Q$ matrices which assigns the orthogonal projection onto $T_zM^m$, we have $\pi^j_{k'}(x):=P_{\vec{\Psi}^j_{k'}(x)}$ and we have
\be
\label{A.I.16}
\pi^j_{k'}\longrightarrow P_{\vec{\Psi}_\infty}\quad\quad\mbox{ strongly in }W^{1,2}_{loc}(S^j_\infty\setminus \cup_{l=1}^{N^j} B_\delta(a^{j,l}))
\ee
  On $S^j_\infty\setminus \cup_{l=1}^{N^j} B_\delta(a^{j,l})$  we denote  $\vec{u}_{i,k'}^{\,\delta}(x):=\pi^j_{k}(x)(\vec{w}_{i,k}^{\,\delta})$. Because of (\ref{A.I.16}) we have
 \be
\label{A.I.17}
  \vec{u}_{i,k}^{\,\delta}\longrightarrow \vec{w}_i^{\,\delta}\quad\quad\mbox{ strongly in }W^{1,2}(S^j_\infty)\quad.
\ee
Consider now the symmetric matrix
  \[
\begin{array}{l}
\ds D^2 \mbox{Area}(\vec{\Phi}_{k})(\vec{u}_{i,k}^{\,\delta},\vec{u}_{i',k}^{\,\delta})=\\[3mm]
\ds\quad\quad\sum_{j=1}^{\mbox{card}(J)}\int_{S^j_\infty} \lf[\lf<d\vec{u}_{i,k}^{\,\delta}\,;\,d\vec{u}_{i',k}^{\,\delta}\rg>_{g_{\vec{\Psi}^j_{k}}}+\lf<d\vec{\Psi}^j_{k}\,;\,d\vec{u}^{\,\delta}_{i,k}\rg>_{g_{\vec{\Psi}_{k}^j}}\lf<d\vec{\Psi}_{k}^j\,;\,d\vec{u}_{i',k}^{\,\delta}\rg>_{g_{\vec{\Psi}_{k}^j}} \rg]\ d\mbox{vol}_{g_{\vec{\Psi}_{k}^j}}\\[3mm]
\ds- \ 2^{-1} \sum_{j=1}^{\mbox{card}(J)}\int_{S^j_\infty} \lf<d\vec{\Psi}_{k}^j\dot{\otimes}\, d\vec{u}_{i,k}^{\,\delta}+d\vec{u}_{i,k}^{\,\delta}\dot{\otimes}\, d\vec{\Psi}_{k}^j, d\vec{\Psi}_{k}^j\dot{\otimes}\, d\vec{u}_{i',k}^{\,\delta}+d\vec{u}^{\,\delta}_{i',k}\dot{\otimes}\, d\vec{\Psi}_{k}^j\rg>\ d\mbox{vol}_{g_{\vec{\Psi}_{k}^j}}
\end{array}
\]
Let $f$ and $g$ be two smooth functions supported on $\vec{\Psi}_\infty(S_\infty^j\setminus\cup_{l=1}^{N^j} B_\delta(a^{j,l}) )$ then one has 
\[
\int_{S^j_\infty}<d(f(\vec{\Psi}_k)), d(g(\vec{\Psi}_k))>_{g_{\vec{\Psi}_{k}}}\ d\mbox{vol}_{g_{\vec{\Psi}_{k}}}=\int_{S^j_\infty}<d(f(\vec{\Psi}_k)), d(g(\vec{\Psi}_k))>_{h^j_k}\ d\mbox{vol}_{h^j_k}
\]
And since $h_{k}^j$ converges in any norms towards $h^j_\infty$, because of the strong $W^{1,2}$ convergence of $\vec{\Psi}_k$ on $S_\infty^j\setminus\cup_{l=1}^{N^j} B_\delta(a^{j,l})$ one has
\be
\label{A.I.18}
\int_{S^j_\infty}<d(f(\vec{\Psi}_k)), d(g(\vec{\Psi}_k))>_{g_{\vec{\Psi}_{k}}}\ d\mbox{vol}_{g_{\vec{\Psi}_{k}}}\longrightarrow \int_{S^j_\infty}<d(f(\vec{\Psi}_\infty)), d(g(\vec{\Psi}_\infty))>_{g_{\vec{\Psi}_{\infty}}}\ d\mbox{vol}_{g_{\vec{\Psi}_{\infty}}}
\ee
In a conformal chart  for $h_{k}^j$ we denote $e^{\la^j_{k'}}:=|\p_{x_1}\vec{\Psi}_{k'}^j|=|\p_{x_2}\vec{\Psi}_{k'}^j|$. Because of the strong $W^{1,2}$ convergence (\ref{A.I.15}) we have
\[
e^{\la^j_{k'}}\longrightarrow e^{\la^j_\infty}=|\p_{x_1}\vec{\Psi}_\infty|=|\p_{x_2}\vec{\Psi}_\infty|\quad\mbox{ a. e. in }\quad S^j_\infty\quad.
\]
Since $e^{\la^j_\infty}>0$ almost everywhere on $S^j_\infty$ we have $e^{-\la^j_{k'}}\longrightarrow e^{-\la^j_\infty}$ almost everywhere and then for $i=1,2$
\[
\p_{x_i}\vec{\Psi}^j_{k}/e^{\la^j_{k}}\longrightarrow \p_{x_i}\vec{\Psi}^j_{\infty}/e^{\la^j_{\infty}}\quad\mbox{ almost everywhere}
\]
Let $f$ ,$g$, $\phi$ and $\psi$ be 4 smooth functions on $M^m$ where $f$ and $g$ are supported on  $\vec{\Psi}_\infty(S_\infty^j\setminus\cup_{l=1}^{N^j} B_\delta(a^{j,l}) )$  one has in local conformal coordinates
\[
\begin{array}{l}
\ds<d(f(\vec{\Psi}^j_k))\otimes d(\phi(\vec{\Psi}^j_k)), d(g(\vec{\Psi}^j_k))\otimes d(\psi(\vec{\Psi}^j_k))>_{g_{\vec{\Psi}^j_{k}}}\ d\mbox{vol}_{g_{\vec{\Psi}^j_{k}}}=\\[3mm]
\ds\sum_{\mu,\nu=1,2} e^{-2\la^j_k} \p_{x_\mu}f(\vec{\Psi}^j_k)\ \p_{x_\nu}\phi(\vec{\Psi}^j_k)\ \p_{x_\mu}g(\vec{\Psi}^j_k)\ \p_{x_\nu}\psi(\vec{\Psi}^j_k)\ dx_1\wedge dx_2
\end{array}
\]
Because of the above
\[
e^{-2\la^j_k} \p_{x_\mu}f(\vec{\Psi}^j_k)\ \p_{x_\nu}\phi(\vec{\Psi}^j_k)\ \p_{x_\mu}g(\vec{\Psi}^j_k)\ \p_{x_\nu}\psi(\vec{\Psi}^j_k)\longrightarrow 
e^{-2\la^j_\infty} \p_{x_\mu}f(\vec{\Psi}^j_\infty)\ \p_{x_\nu}\phi(\vec{\Psi}^j_\infty)\ \p_{x_\mu}g(\vec{\Psi}^j_\infty)\ \p_{x_\nu}\psi(\vec{\Psi}^j_\infty)
\]
almost everywhere and we have moreover
\[
|e^{-2\la^j_k} \p_{x_\mu}f(\vec{\Psi}^j_k)\ \p_{x_\nu}\phi(\vec{\Psi}^j_k)\ \p_{x_\mu}g(\vec{\Psi}^j_k)\ \p_{x_\nu}\psi(\vec{\Psi}^j_k)|\le C\ |\nabla\vec{\Psi}^j_k|^2\rightarrow |\nabla\vec{\Psi}^j_\infty|^2\quad\mbox{strongly in }L^1
\]
Hence the generalized dominated convergence theorem implies
\[
\begin{array}{l}
\ds\int_{S^j_\infty}<d(f(\vec{\Psi}^j_k))\otimes d(\phi(\vec{\Psi}^j_k)), d(g(\vec{\Psi}^j_k))\otimes d(\psi(\vec{\Psi}^j_k))>_{g_{\vec{\Psi}^j_{k}}}\ d\mbox{vol}_{g_{\vec{\Psi}^j_{k}}}\\[3mm]
\ds\longrightarrow\quad\int_{S^j_\infty}<d(f(\vec{\Psi}^j_\infty))\otimes d(\phi(\vec{\Psi}^j_\infty)), d(g(\vec{\Psi}^j_\infty))\otimes d(\psi(\vec{\Psi}^j_\infty))>_{g_{\vec{\Psi}^j_{\infty}}}\ d\mbox{vol}_{g_{\vec{\Psi}^j_{\infty}}}
\end{array}
\]
Similarly we also have
\[
\begin{array}{l}
\ds\int_{S^j_\infty}\lf<d(f(\vec{\Psi}^j_k)), d(g(\vec{\Psi}^j_k))\rg>_{g_{\vec{\Psi}_k}}\ \lf<d(\phi(\vec{\Psi}^j_k)), d(\psi(\vec{\Psi}^j_k))\rg>_{g_{\vec{\Psi}_k}}  \ d\mbox{vol}_{g_{\vec{\Psi}^j_{k}}}\\[3mm]
\ds\longrightarrow\quad\int_{S^j_\infty}\lf<d(f(\vec{\Psi}^j_\infty)), d(g(\vec{\Psi}^j_\infty))\rg>_{g_{\vec{\Psi}_\infty}}\ \lf<d(\phi(\vec{\Psi}^j_\infty)), d(\psi(\vec{\Psi}^j_\infty))\rg>_{g_{\vec{\Psi}_\infty}}  \ d\mbox{vol}_{g_{\vec{\Psi}^j_{\infty}}}
\end{array}
\]
We have
\[
\vec{u}_{i,k}:=\sum_{s=1}^N\chi_s(\vec{\Psi}_k(x))\sum_{t=1}^{n_{s}}P_{\vec{\Psi}_k}\lf(\vec{v}_{i,s}^{\,t}(\vec{\Psi}_k(x))\rg)\ \om_s^t
\]
Because of (\ref{canc}) we have obviously from (\ref{deuxi-18}) that for any choice of $s,t,s',t'$
\[
\begin{array}{l}
\ds D^2 F(\vec{\Phi})\lf(\chi_s(\vec{\Psi}_k(x))P_{\vec{\Psi}_k}\lf(\vec{v}_{i,s}^{\,t}(\vec{\Psi}_k(x))\rg)\ \ti{\om}_s^t, \chi_{s'}(\vec{\Psi}_k(x))P_{\vec{\Psi}_k}\lf(\vec{v}_{i,s'}^{\,t'}(\vec{\Psi}_k(x))\rg)\ \ti{\om}_{s'}^{t'}\rg)\\[5mm]
\ds\quad\le \int_{\ti{\Om}_s^t\cap \ti{\Om}_{s'}^{t'}}(1+|\vec{\mathbb I}_{\vec{\Phi}}|^2_{g_{\vec{\Phi}}})^2\ dvol_{g_{\vec{\Phi}}}
\end{array}
\]
Combining all the above gives
\be
\label{A.I.19}
 D^2 \mbox{Area}(\vec{\Phi}_{k})(\vec{u}_{i,k}^{\,\delta},\vec{u}_{i',k}^{\,\delta})\quad\longrightarrow\quad D^2 \mbox{Area}(\vec{\Phi}_{\infty})(\vec{w}_{i}^{\,\delta},\vec{w}_{i'}^{\,\delta})
\ee
Hence, for $k$ large enough $(D^2 A(\vec{\Phi}_{k})(\vec{u}_{i,k}^{\,\delta},\vec{u}_{i',k}^{\,\delta}))_{i,i'=1\cdots N}$ defines a strictly negative quadratic form.

\medskip

Using now lemma~\ref{deuxi} below we deduce that for any $i,i'\in \{1\cdots N\}$
\be
\label{A.I.20}
\sigma_k^2 \lf|D^2 F(\vec{\Phi}_{k})(\vec{u}_{i,k}^{\,\delta},\vec{u}_{i',k}^{\,\delta})\rg|\le C\sigma_k^2 \lf[ F(\vec{\Phi}_{k})+\mbox{Area}(\vec{\Phi}_{k})^{1/4}\ F(\vec{\Phi}_{k})^{3/4}\rg]=o(1)
\ee
Combining (\ref{A.I.19}) and (\ref{A.I.20}) we obtain that for $k$ large enough $(D^2 A^{\sigma_k}(\vec{\Phi}_{k})(\vec{u}_{i,k}^{\,\delta},\vec{u}_{i',k}^{\,\delta}))_{i,i'=1\cdots N}$ defines a strictly negative quadratic form. This implies inequality (\ref{A.I.14}) and theorem~\ref{th-A.3} is proved.\hfill $\Box$

 \renewcommand{\theequation}{A.\arabic{equation}}
\renewcommand{\theTh}{A.\arabic{Th}}
\renewcommand{\theProp}{A.\arabic{Prop}}
\renewcommand{\theLma}{A.\arabic{Lma}}
\renewcommand{\theCo}{A.\arabic{Co}}
\renewcommand{\theRm}{A.\arabic{Rm}}
\renewcommand{\theequation}{A.\arabic{equation}}
\setcounter{equation}{0} 
\reset
\appendix
\section{Appendix}

\begin{Lma}
\label{deuxi}
Let $M^m$ be a closed sub-manifold of the euclidian space ${\R}^Q$. For any $W^{2,4}-$immersion $\vec{\Phi}$ of an oriented closed surface $\Sigma$ we denote
\[
F(\vec{\Phi}):=\int_\Sigma\lf(1+|\vec{\mathbb I}_{\vec{\Phi}}|^2_{g_{\vec{\Phi}}}\rg)^2\ dvol_{g_{\vec{\Phi}}}
\]
where $\vec{\mathbb I}_{\vec{\Phi}}$ is the second fundamental for of the immersion into $M^m$. The lagrangian $F$ is $C^2$ and there exists a constant $C$ depending only on $M^m$
such that for any perturbation $\vec{w}$ of the form $\vec{v}\circ\vec{\Phi}$ one has
\be
\label{deuxi-der1}
|DF(\vec{\Phi})(\vec{v}(\vec{\Phi}))|\le C\,\int_{\Sigma}\lf(1+|\vec{\mathbb I}_{\vec{\Phi}}|^2_{g_{\vec{\Phi}}}\rg)\ \lf[ (1+|\vec{\mathbb I}_{\vec{\Phi}}|^2_{g_{\vec{\Phi}}})\ |\p \vec{v}|(\vec{\Phi})+ |\vec{\mathbb I}_{\vec{\Phi}}|_{g_{\vec{\Phi}}}\ |\p^2 \vec{v}|(\vec{\Phi})\rg]\ dvol_{g_{\vec{\Phi}}}
\ee
and
\be
\label{deuxi-der2}
|D^2F(\vec{\Phi})(\vec{v}(\vec{\Phi}),\vec{v}(\vec{\Phi}))|\le C\,\int_{\Sigma}\lf(1+|\vec{\mathbb I}_{\vec{\Phi}}|^2_{g_{\vec{\Phi}}}\rg)\ \lf[ (1+|\vec{\mathbb I}_{\vec{\Phi}}|^2_{g_{\vec{\Phi}}})\ |\p \vec{v}|^2(\vec{\Phi})+\ |\p^2 \vec{v}|^2(\vec{\Phi})\rg]\ dvol_{g_{\vec{\Phi}}}
\ee
\hfill $\Box$
\end{Lma}
\noindent{\bf Proof of lemma~\ref{deuxi}.}

In local coordinates we  denote the second fundamental form
\[
\vec{\mathbb I}_{\vec{\Phi}}=\pi_{\vec{n}}\lf(d^2\vec{\Phi}   \rg)=\pi_{\vec{n}}\lf(\p^2_{x_ix_j}\vec{\Phi}\rg)\ dx_i\otimes dx_j
\]
we have
\be
\label{deuxi-1}
|\vec{\mathbb I}_{\vec{\Phi}}|^2_{g_{\vec{\Phi}}}:=\lf|\pi_{\vec{n}}\lf(d^2\vec{\Phi}   \rg)\rg|^2_{g_{\vec{\Phi}}}=\sum_{i,j,k,l}g^{ik}g^{jl}\pi_{\vec{n}}\p^2_{x_ix_j}\vec{\Phi}\cdot\pi_{\vec{n}}\p^2_{x_kx_l}\vec{\Phi}
\ee
Denote $\pi_T$ the projection onto the tangent plane to the immersion. We have in local coordinates
\be
\label{deuxi-2}
\pi_T(\vec{X})=\sum_{i,j=1}^2 g^{ij}\ \p_{x_i}\vec{\Phi}\cdot\vec{X}\ \p_{x_j}\vec{\Phi}
\ee
Hence
\be
\label{deuxi-3}
\lf.\pi_{\vec{n}}\frac{d\pi_{\vec{n}}}{dt}\rg|_{t=0}(\vec{X})=-\sum_{i,j=1}^2 g^{ij}\ \p_{x_i}\vec{\Phi}\cdot\vec{X}\ \pi_{\vec{n}}\lf(\p_{x_j}\vec{w}\rg)
\ee
We have clearly
\[
\frac{d g_{ij}}{dt}=\p_{x_i}\vec{\Phi}\cdot\p_{x_j}\vec{w}+\p_{x_i}\vec{w}\cdot\p_{x_j}\vec{\Phi}
\]
Hence
\be
\label{deuxi-4}
\frac{d g^{ij}}{dt}=-\,g^{ik} g^{jl}\,\lf[\p_{x_k}\vec{\Phi}\cdot\p_{x_l}\vec{w}+\p_{x_k}\vec{w}\cdot\p_{x_l}\vec{\Phi}\rg]:= -2\,(d\vec{\Phi}\dot{\otimes}_Sd\vec{w})^{ij}
\ee
We have then
\be
\label{deuxi-5}
\lf.\frac{d|\vec{\mathbb I}_{\vec{\Phi}}|^2_{g_{\vec{\Phi}}}}{dt}\rg|_{t=0}=2\, \lf<\pi_{\vec{n}}\lf(d^2\vec{\Phi}   \rg),\pi_{\vec{n}}\lf(D^{g_{\vec{\Phi}}}d\vec{w}\rg)\rg>_{g_{\vec{\Phi}}}-\,4\lf( g\otimes(d\vec{\Phi}\dot{\otimes}_Sd\vec{w})\res \vec{{\mathbb I}}_{\vec{\Phi}}\dot{\otimes}\vec{{\mathbb I}}_{\vec{\Phi}}\rg)
\ee
where $\res$ is the contraction operator between  $4-$contravariant and $4-$covariant tensors and
\be
\label{deuxi-6}
D^{g_{\vec{\Phi}}}d\vec{w}:=\lf[\p^2_{x_ix_j}\vec{w} -\sum_{rs=1}^2g^{rs}\p_{x_r}\vec{\Phi}\cdot \p^2_{x_ix_j}\vec{\Phi}\ \p_{x_s}\vec{w}\rg]\, dx_i\otimes dx_j\quad.
\ee
This gives in particular that
\be
\label{deuxi-7}
\begin{array}{l}
\ds\lf.\frac{d}{dt}\int_{\Sigma}(1+|\vec{\mathbb I}_{\vec{\Phi}}|^2_{g_{\vec{\Phi}}})^2\ dvol_{g_{\vec{\Phi}}}\rg|_{t=0}=DF(\vec{\Phi})(\vec{w})\\[5mm]
\ds\quad=4\, \int_{\Sigma}(1+|\vec{\mathbb I}_{\vec{\Phi}}|^2_{g_{\vec{\Phi}}}) \lf[\lf<\vec{\mathbb I}_{\vec{\Phi}},D^{g_{\vec{\Phi}}}d\vec{w}\rg>_{g_{\vec{\Phi}}}-\,2\lf( g\otimes(d\vec{\Phi}\dot{\otimes}_Sd\vec{w})\rg)\res \lf(\vec{{\mathbb I}}_{\vec{\Phi}}\dot{\otimes}\vec{{\mathbb I}}_{\vec{\Phi}}\rg)\rg]\ dvol_{g_{\vec{\Phi}}}\\[5mm]
\ds\quad+\int_{\Sigma}(1+|\vec{\mathbb I}_{\vec{\Phi}}|^2_{g_{\vec{\Phi}}})^2 \lf<d\vec{\Phi};d\vec{w}\rg>_{g_{\vec{\Phi}}}\ dvol_{g_{\vec{\Phi}}}
\end{array}
\ee
For $\vec{w}:=\vec{v}\lf(\vec{\Phi}\rg)$ we have 
\[
\begin{array}{l}
\ds\p^2_{x_ix_j}\vec{w} -\sum_{rs=1}^2g^{rs}\p_{x_r}\vec{\Phi}\cdot \p^2_{x_ix_j}\vec{\Phi}\ \p_{x_s}\vec{w}=\sum_{\al,\beta=1}^Q\p^2_{z_\al z_\beta}\vec{v}(\vec{\Phi})\, \p_{x_i}\vec{\Phi}^{\,\al}\,\p_{x_j}\vec{\Phi}^{\,\beta}\\[5mm]
\ds+\sum_{\al=1}^Q\p_{z_\al }\vec{v}(\vec{\Phi})\, \lf[\p^2_{x_ix_j}\vec{\Phi}^{\,\al}
-\ds\sum_{rs=1}^2g^{rs}\p_{x_r}\vec{\Phi}\cdot \p^2_{x_ix_j}\vec{\Phi}\  \p_{x_s}\vec{\Phi}^{\,\al}\rg]
\end{array}
\]
We have
\be
\label{deuxi-8}
\pi_{T}\lf(\p^2_{x_ix_j}\vec{\Phi}
\rg)=\sum_{rs=1}^2g^{rs}\,\p^2_{x_ix_j}\vec{\Phi}\cdot\p_{x_r}\vec{\Phi}\  \p_{x_s}\vec{\Phi}\quad.
\ee
Hence
\[
\p^2_{x_ix_j}\vec{\Phi}-\sum_{rs=1}^2g^{rs}\p_{x_r}\vec{\Phi}\cdot \p^2_{x_ix_j}\vec{\Phi}\  \p_{x_s}\vec{\Phi}=\pi_{\vec{n}}(\p^2_{x_ix_j}\vec{\Phi})=\vec{\mathbb I}_{ij}\quad.
\]
This implies that
\be
\label{deuxi-9}
D^{g_{\vec{\Phi}}}d\vec{w}=\sum_{\al,\beta=1}^Q\p^2_{z_\al z_\beta}\vec{v}(\vec{\Phi})\, d\vec{\Phi}^{\,\al}\otimes d \vec{\Phi}^{\,\beta}+\sum_{\al=1}^Q\p_{z_\al }\vec{v}(\vec{\Phi})\, \vec{\mathbb I}_{ij}^{\, \al}
\ee
We deduce
\be
\label{deuxi-10}
|DF(\vec{\Phi})(\vec{v}(\vec{\Phi}))|\le C\,\int_{\Sigma}\lf(1+|\vec{\mathbb I}_{\vec{\Phi}}|^2_{g_{\vec{\Phi}}}\rg)\ \lf[ (1+|\vec{\mathbb I}_{\vec{\Phi}}|^2_{g_{\vec{\Phi}}})\ |\p \vec{v}|(\vec{\Phi})+ |\vec{\mathbb I}_{\vec{\Phi}}|_{g_{\vec{\Phi}}}\ |\p^2 \vec{v}|(\vec{\Phi})\rg]\ dvol_{g_{\vec{\Phi}}}
\ee
We now compute the second derivative
\be
\label{deuxi-11}
\begin{array}{l}
\ds\lf.\frac{d}{dt}DF(\vec{\Phi}_t)(\vec{w})\rg|_{t=0}\\[5mm]
\ds=4\, \int_{\Sigma}(1+|\vec{\mathbb I}_{\vec{\Phi}}|^2_{g_{\vec{\Phi}}}) \lf[\lf<\vec{\mathbb I}_{\vec{\Phi}},D^{g_{\vec{\Phi}}}d\vec{w}\rg>_{g_{\vec{\Phi}}}-\,2\lf( g\otimes(d\vec{\Phi}\dot{\otimes}_Sd\vec{w})\rg)\res \lf(\vec{{\mathbb I}}_{\vec{\Phi}}\dot{\otimes}\vec{{\mathbb I}}_{\vec{\Phi}}\rg)\rg]\ \lf<d\vec{\Phi};d\vec{w}\rg>_{g_{\vec{\Phi}}}\ dvol_{g_{\vec{\Phi}}}\\[5mm]
\ds+\int_{\Sigma}|1+|\vec{\mathbb I}_{\vec{\Phi}}|^2_{g_{\vec{\Phi}}}|^2 \lf|\lf<d\vec{\Phi};d\vec{w}\rg>_{g_{\vec{\Phi}}}\rg|^2+8\, \lf|\lf<\vec{\mathbb I}_{\vec{\Phi}},D^{g_{\vec{\Phi}}}d\vec{w}\rg>_{g_{\vec{\Phi}}}-\,2\lf( g\otimes(d\vec{\Phi}\dot{\otimes}_Sd\vec{w})\rg)\res \lf(\vec{{\mathbb I}}_{\vec{\Phi}}\dot{\otimes}\vec{{\mathbb I}}_{\vec{\Phi}}\rg)\rg|^2\ dvol_{g_{\vec{\Phi}}}\\[5mm]
\ds+\,4\, \int_{\Sigma}(1+|\vec{\mathbb I}_{\vec{\Phi}}|^2_{g_{\vec{\Phi}}}) \frac{d}{dt}\lf[\lf<\vec{\mathbb I}_{\vec{\Phi}_t},D^{g_{\vec{\Phi}_t}}d\vec{w}\rg>_{g_{\vec{\Phi}_t}}-\,2\lf( g_{\vec{\Phi}_t}\otimes(d\vec{\Phi}_t\dot{\otimes}_Sd\vec{w})\rg)\res \lf(\vec{{\mathbb I}}_{\vec{\Phi}_t}\dot{\otimes}\vec{{\mathbb I}}_{\vec{\Phi}_t}\rg)\rg]\ dvol_{g_{\vec{\Phi}}}\\[5mm]
\ds+\int_{\Sigma}(1+|\vec{\mathbb I}_{\vec{\Phi}}|^2_{g_{\vec{\Phi}}})^2 \lf<d\vec{w};d\vec{w}\rg>_{g_{\vec{\Phi}}}\ dvol_{g_{\vec{\Phi}}}-2\,\int_{\Sigma} (1+|\vec{\mathbb I}_{\vec{\Phi}}|^2_{g_{\vec{\Phi}}})^2 \ \lf(d\vec{\Phi}\dot{\otimes}_S d\vec{w}\rg)\res \lf(d\vec{\Phi}\otimes d\vec{w}\rg)\ dvol_{g_{\vec{\Phi}}}
\end{array}
\ee
We have in one hand
\be
\label{deuxi-12}
\pi_{\vec{n}}\frac{d}{dt}\vec{\mathbb I}_{\vec{\Phi}_t}=\pi_{\vec{n}}\lf(D^{g_{\vec{\Phi}}}d\vec{w}\rg)
\ee
in the other hand
\be
\label{deuxi-13}
\begin{array}{l}
\ds\frac{d}{dt}\lf( \sum_{r=1}^2g^{rs}_{\vec{\Phi}_t}\ \p_{x_r}\vec{\Phi}_t\cdot \p^2_{x_ix_j}\vec{\Phi}_t \rg)=\sum_{r=1}^2g^{rs}\ \p_{x_r}\vec{\Phi}\cdot \p^2_{x_ix_j}\vec{w}+g^{rs}\ \p_{x_r}\vec{w}\cdot \p^2_{x_ix_j}\vec{\Phi}+\frac{d g^{rs}}{dt} \p_{x_r}\vec{\Phi}\cdot \p^2_{x_ix_j}\vec{\Phi}\\[5mm]
\ds\quad=\sum_{r=1}^2g^{rs}\ \p_{x_r}\vec{\Phi}\cdot \p^2_{x_ix_j}\vec{w}+g^{rs}\ \p_{x_r}\vec{w}\cdot \pi_{\vec{n}}(\p^2_{x_ix_j}\vec{\Phi})+\frac{d g^{rs}}{dt} \p_{x_r}\vec{\Phi}\cdot \p^2_{x_ix_j}\vec{\Phi}\\[5mm]
\ds\quad+\sum_{r,k,l=1}^2g^{rs}\ g^{kl}\ \p_{x_r}\vec{w}\cdot \p_{x_k}\vec{\Phi}\ \p_{x_l}\vec{\Phi}\cdot\p^2_{x_ix_j}\vec{\Phi}
\end{array}
\ee
we have
\be
\label{deuxi-14}
\begin{array}{l}
\ds\sum_{r=1}^2\frac{d g^{rs}}{dt} \p_{x_r}\vec{\Phi}\cdot \p^2_{x_ix_j}\vec{\Phi}=- \sum_{r,k,l=1}^2 g^{rk}\,g^{sl}\,\p_{x_r}\vec{\Phi}\cdot \p^2_{x_ix_j}\vec{\Phi}\, \lf[\p_{x_k}\vec{w}\cdot\p_{x_l}\vec{\Phi}+ \p_{x_l}\vec{w}\cdot\p_{x_k}\vec{\Phi} \rg]\\[5mm]
\ds=- \sum_{r,k,l=1}^2 g^{lk}\,g^{sr}\,\p_{x_l}\vec{\Phi}\cdot \p^2_{x_ix_j}\vec{\Phi}\ \p_{x_k}\vec{w}\cdot\p_{x_r}\vec{\Phi}     +  g^{lk}\,g^{sr}\,\p_{x_l}\vec{\Phi}\cdot \p^2_{x_ix_j}\vec{\Phi}\ \p_{x_r}\vec{w}\cdot\p_{x_k}\vec{\Phi}
\end{array}
\ee
Combining (\ref{deuxi-13}) and (\ref{deuxi-14}) we obtain
\be
\label{deuxi-15}
\ds\frac{d}{dt}\lf( \sum_{r=1}^2g^{rs}_{\vec{\Phi}_t}\ \p_{x_r}\vec{\Phi}_t\cdot \p^2_{x_ix_j}\vec{\Phi}_t \rg)=\sum_{r=1}^2g^{rs}\ \p_{x_r}\vec{\Phi}\cdot (D^{g_{\vec{\Phi}}}d\vec{w})_{ij}+\sum_{r=1}^2 g^{rs}\, \p_{x_r}\vec{w}\cdot\vec{\mathbb I}_{ij}
\ee
Thus
\be
\label{deuxi-16}
\ds\frac{d}{dt}\lf(D^{g_{\vec{\Phi}_t}}d\vec{w}\rg)=-\sum_{i,j=1}^2\lf[ \sum_{r=1}^2g^{rs}\ \p_{x_r}\vec{\Phi}\cdot (D^{g_{\vec{\Phi}}}d\vec{w})_{ij}\ \p_{x_s}\vec{w}+ g^{rs}\, \p_{x_r}\vec{w}\cdot\vec{\mathbb I}_{ij}\ \p_{x_s}\vec{w}\rg] \ dx_i\otimes dx_j
\ee
Combining (\ref{deuxi-4}), (\ref{deuxi-12}) and (\ref{deuxi-16}) we obtain
\be
\label{deuxi-17}
\begin{array}{l}
\ds \frac{d}{dt}\lf[\lf<\vec{\mathbb I}_{\vec{\Phi}_t},D^{g_{\vec{\Phi}_t}}d\vec{w}\rg>_{g_{\vec{\Phi}_t}}-\,2\lf( g_{\vec{\Phi}_t}\otimes(d\vec{\Phi}_t\dot{\otimes}_Sd\vec{w})\rg)\res \lf(\vec{{\mathbb I}}_{\vec{\Phi}_t}\dot{\otimes}\vec{{\mathbb I}}_{\vec{\Phi}_t}\rg)\rg]\\[5mm]
\ds=\lf|\pi_{\vec{n}}\lf(D^gd\vec{w}\rg)\rg|_{g_{\vec{\Phi}}}^2-\lf<\vec{\mathbb I}\ ;\sum_{i,j=1}^2g^{ij}\p_{x_i}\vec{\Phi}\cdot D^g d\vec{w}\ \p_{x_j}\vec{w}\rg>_{g_{\vec{\Phi}}}+ 4\, 
\lf[(d\vec{\Phi}\otimes_S d\vec{w})\otimes(d\vec{\Phi}\otimes_S d\vec{w})\rg]\res (\vec{\mathbb I}\dot{\otimes}\vec{\mathbb I})\\[5mm]
\ds-4\, \lf(g\otimes (d\vec{\Phi}\otimes_S d\vec{w})\rg)\res \lf(\vec{\mathbb I}\otimes \pi_{\vec{n}}(D^g d\vec{w})+\pi_{\vec{n}}(D^g d\vec{w})\otimes\vec{\mathbb I}\rg)-2\  \lf(g\otimes (d\vec{w}\otimes_S d\vec{w})\rg)\res (\vec{\mathbb I}\dot{\otimes}\vec{\mathbb I}) 
\end{array}
\ee
Combining (\ref{deuxi-11}) and (\ref{deuxi-17}) gives 
\be
\label{deuxi-18}
\begin{array}{l}
\ds D^2F(\vec{\Phi})(\vec{w},\vec{w})=\\[5mm]
\ds4\, \int_{\Sigma}(1+|\vec{\mathbb I}_{\vec{\Phi}}|^2_{g_{\vec{\Phi}}}) \lf[\lf<\vec{\mathbb I}_{\vec{\Phi}},D^{g_{\vec{\Phi}}}d\vec{w}\rg>_{g_{\vec{\Phi}}}-\,2\lf( g_{\vec{\Phi}}\otimes(d\vec{\Phi}\dot{\otimes}_Sd\vec{w})\rg)\res \lf(\vec{{\mathbb I}}_{\vec{\Phi}}\dot{\otimes}\vec{{\mathbb I}}_{\vec{\Phi}}\rg)\rg]\ \lf<d\vec{\Phi};d\vec{w}\rg>_{g_{\vec{\Phi}}}\ dvol_{g_{\vec{\Phi}}}\\[5mm]
\ds+\int_{\Sigma}|1+|\vec{\mathbb I}_{\vec{\Phi}}|^2_{g_{\vec{\Phi}}}|^2 \lf|\lf<d\vec{\Phi};d\vec{w}\rg>_{g_{\vec{\Phi}}}\rg|^2+8\, \lf|\lf<\vec{\mathbb I}_{\vec{\Phi}},D^{g_{\vec{\Phi}}}d\vec{w}\rg>_{g_{\vec{\Phi}}}-\,2\lf( g\otimes(d\vec{\Phi}\dot{\otimes}_Sd\vec{w})\rg)\res \lf(\vec{{\mathbb I}}_{\vec{\Phi}}\dot{\otimes}\vec{{\mathbb I}}_{\vec{\Phi}}\rg)\rg|^2\ dvol_{g_{\vec{\Phi}}}\\[5mm]
\ds+\,4\, \int_{\Sigma}(1+|\vec{\mathbb I}_{\vec{\Phi}}|^2_{g_{\vec{\Phi}}}) \lf[\lf|\pi_{\vec{n}}\lf(D^{g_{\vec{\Phi}}}d\vec{w}\rg)\rg|_{g_{\vec{\Phi}}}^2-\lf<\vec{\mathbb I}_{\vec{\Phi}}\ ;\sum_{i,j=1}^2g^{ij}\p_{x_i}\vec{\Phi}\cdot D^g_{g_{\vec{\Phi}}} d\vec{w}\ \p_{x_j}\vec{w}\rg>_{g_{\vec{\Phi}}}\rg] \ dvol_{g_{\vec{\Phi}}}\\[5mm]
\ds+\,16\, \int_{\Sigma}(1+|\vec{\mathbb I}_{\vec{\Phi}}|^2_{g_{\vec{\Phi}}})   \lf[(d\vec{\Phi}\otimes_S d\vec{w})\otimes(d\vec{\Phi}\otimes_S d\vec{w})\rg]\res (\vec{\mathbb I}_{\vec{\Phi}}\dot{\otimes}\vec{\mathbb I}_{\vec{\Phi}})          \ dvol_{g_{\vec{\Phi}}}\\[5mm]
\ds-\,16\, \int_{\Sigma}(1+|\vec{\mathbb I}_{\vec{\Phi}}|^2_{g_{\vec{\Phi}}})  \  \lf(g_{\vec{\Phi}}\otimes (d\vec{\Phi}\otimes_S d\vec{w})\rg)\res \lf(\vec{\mathbb I}_{\vec{\Phi}}\otimes \pi_{\vec{n}}(D^{g_{\vec{\Phi}}} d\vec{w}) +\pi_{\vec{n}}(D^{g_{\vec{\Phi}}} d\vec{w})\otimes\vec{\mathbb I}_{\vec{\Phi}}\rg)  \ dvol_{g_{\vec{\Phi}}}\\[5mm]
\ds-\, 8\, \int_{\Sigma}(1+|\vec{\mathbb I}_{\vec{\Phi}}|^2_{g_{\vec{\Phi}}})  \lf(g_{\vec{\Phi}}\otimes (d\vec{w}\otimes_S d\vec{w})\rg)\res (\vec{\mathbb I}_{\vec{\Phi}}\dot{\otimes}\vec{\mathbb I}_{\vec{\Phi}})  \ dvol_{g_{\vec{\Phi}}}\\[5mm]
\ds+\int_{\Sigma}(1+|\vec{\mathbb I}_{\vec{\Phi}}|^2_{g_{\vec{\Phi}}})^2 \lf<d\vec{w};d\vec{w}\rg>_{g_{\vec{\Phi}}}\ dvol_{g_{\vec{\Phi}}}-2\,\int_{\Sigma} (1+|\vec{\mathbb I}_{\vec{\Phi}}|^2_{g_{\vec{\Phi}}})^2 \ \lf(d\vec{\Phi}\dot{\otimes}_S d\vec{w}\rg)\res \lf(d\vec{\Phi}\otimes d\vec{w}\rg)\ dvol_{g_{\vec{\Phi}}}
\end{array}
\ee
For $\vec{w}:=\vec{v}(\vec{\Phi})$, using (\ref{deuxi-9}), we deduce
\be
\label{deuxi-19}
|D^2F(\vec{\Phi})(\vec{v}(\vec{\Phi}),\vec{v}(\vec{\Phi}))|\le C\,\int_{\Sigma}\lf(1+|\vec{\mathbb I}_{\vec{\Phi}}|^2_{g_{\vec{\Phi}}}\rg)\ \lf[ (1+|\vec{\mathbb I}_{\vec{\Phi}}|^2_{g_{\vec{\Phi}}})\ |\p \vec{v}|^2(\vec{\Phi})+\ |\p^2 \vec{v}|^2(\vec{\Phi})\rg]\ dvol_{g_{\vec{\Phi}}}
\ee
This concludes the proof of lemma~\ref{deuxi}.\hfill $\Box$

\end{document}